
\def\ignore#1{}
 

\newcount\sectnum
\newcount\subsectnum
\newcount\eqnumber

\global\eqnumber=1\sectnum=0


\def\lab{(\the\sectnum.\the\eqnumber)}



\def\show#1{#1}



\def\smskip{\vskip 5 pt}
\def\medskip{\vskip 10 pt}
\def\bigskip{\vskip 15 pt}
\def\pn{\par\noindent}
\def\br{\break}

\def\bl{\bigl} 
\def\br{\bigr} 
\def\lf{\left}
\def\ri{\right}

\def\ol#1{\overline{#1}}

\def\b{\beta}

\def\g{\gamma}
\def\m{\mu}

\def\e{\epsilon}

\def\tl{\tilde}

\def\old#1{}
\def\leaderfill{\leaders\hbox to 1em{\hss.\hss}\hfill}


\parindent=2pc
\baselineskip=15pt
\vsize=8.7 true in
\voffset=0.125 true in
\parskip=3pt


\def\minprob#1#2#3{$$\eqalign{&\hbox{minimize\ \ }#1\cr &\hbox{subject to\ \
}#2\cr}\ifnum 0=#3{}\else\eqno(#3)\fi$$}        
     
\def\maxprob#1#2#3{$$\eqalign{&\hbox{maximize\ \ }#1\cr &\hbox{subject to\ \
}#2\cr}\ifnum 0=#3{}\else\eqno(#3)\fi$$}        
     
\def\aligntwo#1#2#3#4#5{$$\eqalign{#1&#2\cr #3&#4\cr}
\ifnum 0=#5{}\else\eqno(#5)\fi$$}
\def\alignthree#1#2#3#4#5#6#7{$$\eqalign{#1&#2\cr #3&#4\cr #5&#6\cr}
\ifnum 0=#7{}\else\eqno(#7)\fi$$}


\def\eqnum{\eqno{\hbox{(\the\sectnum.\the\eqnumber)}\global\advance\eqnumber
by1}}

\def\eqnu{\eqno{\hbox{(\the\sectnum.\the\eqnumber)}\global\advance\eqnumber
by1}}

\newcount\examplnumber
\def\examplnum{\global\advance\examplnumber by1}

\newcount\figrnumber
\def\figrnum{\global\advance\figrnumber by1}

\newcount\propnumber
\def\propnum{\global\advance\propnumber by1}

\newcount\defnumber
\def\defnum{\global\advance\defnumber by1}

\newcount\lemmanumber
\def\lemmanum{\global\advance\lemmanumber by1}

\newcount\assumptionnumber
\def\assumptionnum{\global\advance\assumptionnumber by1}

\def\exampl{\the\sectnum.\the\examplnumber}
\def\figr{\the\sectnum.\the\figrnumber}
\def\propn{\the\sectnum.\the\propnumber}
\def\defn{\the\sectnum.\the\defnumber}
\def\lemman{\the\sectnum.\the\lemmanumber}
\def\assumptionn{\the\sectnum.\the\assumptionnumber}

\def\section#1{\goodbreak\vskip 3pc plus 6pt minus 3pt\leftskip=-2pc
   \global\advance\sectnum by 1\eqnumber=1
\global\examplnumber=1\figrnumber=1\propnumber=1\defnumber=1\lemmanumber=1\assumptionnumber=1%
   \line{\hfuzz=1pc{\hbox to 3pc{\bf 
   \vtop{\hfuzz=1pc\hsize=38pc\hyphenpenalty=10000\noindent\uppercase{\the\sectnum.\quad #1}}\hss}}
			\hfill}
			\leftskip=0pc\nobreak\tenf
			\vskip 1pc plus 4pt minus 2pt\noindent\ignorespaces}



\def\sect#1{\noindent\leftskip=-2pc\tenf
   \goodbreak\vskip 1pc plus 4pt minus 2pt
                \global\advance\subsectnum by 1\eqnumber=1
   \line{\hfuzz=1pc{\hbox to 3pc{\bf 
   \vtop{\hfuzz=1pc\hsize=38pc\hyphenpenalty=10000\noindent\uppercase{{\bf #1}}}\hss}}
                        \hfill}
   \leftskip=0pc\nobreak\tenf
                        \vskip 1pc plus 4pt minus 2pt\nobreak\noindent\ignorespaces}

\def\subsection#1{\noindent\leftskip=0pc\tenf
   \goodbreak\vskip 1pc plus 4pt minus 2pt
   \line{\hfuzz=1pc{\hbox to 3pc{\bf 
   \vtop{\hfuzz=1pc\hsize=38pc\hyphenpenalty=10000\noindent{\bf #1}}\hss}}
                        \hfill}
   \leftskip=0pc\nobreak\tenf
                        \vskip 1pc plus 4pt minus 2pt\nobreak\noindent\ignorespaces}
\def\subsubsection#1{\goodbreak\vskip 1pc plus 4pt minus 2pt
   \hfuzz=3pc\leftskip=0pc\noindent\tenit #1 \nobreak\tenf\vskip 6pt plus 1pt
                                minus 1pt\nobreak\ignorespaces\leftskip=0pc}
%

\def\beginexample#1{\noindent\goodbreak\vskip 6pt plus 1pt minus 1pt
\noindent
  \hbox {\bf Example #1\hss}
  \nobreak\vskip 4pt plus 1pt minus 1pt \nobreak\noindent\ninef
  \global\advance
                \leftskip by\parindent\pn}
\def\endexample{\vskip 12pt\tenf\par
  \global\advance\leftskip by -\parindent
  }

\def\beginexercise#1{\noindent\goodbreak\vskip 6pt plus 1pt minus 1pt \noindent\global\normalbaselineskip=12pt
  \hbox {\bf Exercise #1\hss}
  \nobreak\vskip 4pt plus 1pt minus 1pt 
  \nobreak\noindent\ninef\global\advance\leftskip
                        by\parindent\pn}
\def\endexercise{\vskip 12pt\tenf\par
  \global\advance\leftskip by -\parindent
  }

\def\beginsection#1{\noindent\goodbreak\vskip 6pt plus 1pt minus 1pt \noindent\global\normalbaselineskip=12pt
  \hbox {\it #1\hss}
  \vskip 0.1pt plus 1pt minus 1pt \nobreak\noindent\ninef\global\advance
                \leftskip by\parindent\noindent\pn}
\def\endsection{\vskip 12pt\tenf\par
  \global\advance\leftskip by -\parindent
}

%


\def\proposition#1{\smskip\pn{\bf Proposition #1}\quad}
\def\proof{\smskip\pn{\bf Proof:}\quad} 
\def\definition#1{\smskip\pn{\bf
Definition #1}\quad}

 \def\qed{\quad{\bf
Q.E.D.} \par\bigskip}
\def\ref{\smskip\pn}

\def\chapter#1#2{{\bf \centerline{\helbigbig
{#1}}}\bigskip\bigskip{\bf \centerline{\helbigbig
{#2}}}\bigskip\bigskip} 



\def\longpapertitle#1#2#3{{\bf \centerline{\helbigb
{#1}}}\bigskip{\bf \centerline{\helbigb
{#2}}}\bigskip\bigskip{\centerline{
by}}\bigskip{\bf \centerline{
{#3}}}\bigskip\bigskip} 


\def\nitem#1{\smskip\item{#1}}

\newcount\alphanum
\newcount\romnum

\def\alphaenumerate{\ifcase\alphanum \or (a)\or (b)\or (c)\or (d)\or (e)\or
(f)\or (g)\or (h)\or (i)\or (j)\or (k)\fi}
\def\romenumerate{\ifcase\romnum \or (i)\or (ii)\or (iii)\or (iv)\or (v)\or
(vi)\or (vii)\or (viii)\or (ix)\or (x)\or (xi)\fi}

\def\alist{\begingroup\vskip10pt\alphanum=1
\parskip=2pt\parindent=0pt \leftskip=3pc
\everypar{\llap{{\rm\alphaenumerate\hskip1em}}\advance\alphanum by1}}

\def\nolist{\begingroup\vskip10pt\alphanum=0
\parskip=2pt\parindent=0pt \leftskip=3pc
\everypar{\llap{\global\advance\alphanum by1(\the\alphanum)\hskip1em}}}

\def\romlist{\begingroup\vskip10pt\romnum=1
\parskip=2pt\parindent=0pt \leftskip=5pc
\everypar{\llap{{\rm\romenumerate\hskip1em}}\advance\romnum by1}}



\long\def\fig#1#2#3{\vbox{\vskip1pc\vskip#1
\prevdepth=12pt \baselineskip=12pt
\vskip1pc
\hbox to\hsize{\hfill\vtop{\hsize=25pc\noindent{\eightbf Figure #2\ }
{\eightpoint#3}}\hfill}}}

\long\def\widefig#1#2#3{\vbox{\vskip1pc\vskip#1
\prevdepth=12pt \baselineskip=12pt
\vskip1pc
\hbox to\hsize{\hfill\vtop{\hsize=28pc\noindent{\eightbf Figure #2\ }
{\eightpoint#3}}\hfill}}}

\long\def\table#1#2{\vbox{\vskip0.5pc
\prevdepth=12pt \baselineskip=12pt
\hbox to\hsize{\hfill\vtop{\hsize=25pc\noindent{\eightbf Table #1\ }
{\eightpoint#2}}\hfill}}}

 
\def\rightheadline#1{\headline{\tenrm\hfil #1}}


\long\def\leftfig#1#2{\vbox{\smskip\hsize=220pt
\vtop{{\noindent {\bf #1}}}
\smskip
\noindent
\vbox{{\noindent #2}}
}}

\long\def\rightfig#1#2#3{\vbox{\smskip\vskip#1
\prevdepth=12pt \baselineskip=12pt
\hsize=210pt
\smskip
\vbox{\noindent{\eightbold #2}
\hskip1em{\eightpoint#3}}
}}

\long\def\concept#1#2#3#4#5{\bigskip\hrule
\vbox{\hbox{\leftfig{#1}{#2} \hskip3em
\rightfig{#3}{#4}{#5}} \smskip}
\hrule\bigskip}


\long\def\bconcept#1#2#3#4#5#6#7{
\vbox{
\hbox to \hsize{\vtop{\par #1}}
\concept{#2}{#3}{#4}{#5}{#6}
\hbox to \hsize{\vtop{\par #7}}
\smskip}
}




\def\boxit#1{\vbox{\hrule\hbox{\vrule\kern3pt
                                \vbox{\kern3pt#1\kern3pt}\kern3pt\vrule}\hrule}}
\def\centerboxit#1{$$\vbox{\hrule\hbox{\vrule\kern3pt
                                \vbox{\kern3pt#1\kern3pt}\kern3pt\vrule}\hrule}$$}

\long\def\boxtext#1#2{$$\boxit{\vbox{\hsize #1\noindent\strut #2\strut}}$$}

%
%
%

\def\picture #1 by #2 (#3){
  \vbox to #2{
    \hrule width #1 height 0pt depth 0pt
    \vfill
    \special{picture #3} 
    }
  }

\def\scaledpicture #1 by #2 (#3 scaled #4){{
  \dimen0=#1 \dimen1=#2
  \divide\dimen0 by 1000 \multiply\dimen0 by #4
  \divide\dimen1 by 1000 \multiply\dimen1 by #4
  \picture \dimen0 by \dimen1 (#3 scaled #4)}
  }

%
%

\long\def\captfig#1#2#3#4#5{\vbox{\vskip1pc
\hbox to\hsize{\hfill{\picture #1 by #2 (#3)}\hfill}
\prevdepth=9pt \baselineskip=9pt
\vskip1pc
\hbox to\hsize{\hfill\vtop{\hsize=24pc\noindent{\eightbold Figure #4}
\hskip1em{\eightpoint#5}}\hfill}}}

%
%
%

\def\illustration #1 by #2 (#3){
  \vskip#2\hskip#1\special{illustration #3} 
    }

\def\scaledillustration #1 by #2 (#3 scaled #4){{
  \dimen0=#1 \dimen1=#2
  \divide\dimen0 by 1000 \multiply\dimen0 by #4
  \divide\dimen1 by 1000 \multiply\dimen1 by #4
  \illustration \dimen0 by \dimen1 (#3 scaled #4)}
  }


\newbox\graybox
\newdimen\xgrayspace
\newdimen\ygrayspace
%
%
%
%
%
%
%
%
%

\def\Textshade#1#2#3#4#5#6{%
    \xgrayspace=#4pt%
    \ygrayspace=#4pt%
    \def\grayshade{#3}%
    \def\linewidth{#5}%
    \def\theradius{#6}%
    \setbox\graybox=\hbox{\surroundboxa{#2}}%
    \hbox{%
    \hbox to 0pt{%
    \PScommands
}%
    \box\graybox}}%
%
%
\long%

\long%
\def\Parashade#1#2#3#4#5#6#7{%
    \xgrayspace=#4pt%
    \ygrayspace=#4pt%
    \def\grayshade{#3}%
    \def\linewidth{#5}%
    \def\theradius{#6}%
    \def\thevskip{#7pt}%
    \setbox\graybox=\hbox{\surroundboxb{#2}}%
    \vskip\thevskip%
    \hbox{%
    \hbox to 0pt{%
    \PScommands
}%
     \box\graybox}%
     \vskip\thevskip%
}%
%
%
%
\long\def\surroundboxa#1{\leavevmode\hbox{\vtop{%
\vbox{\kern\ygrayspace%
\hbox{\kern\xgrayspace#1%
      \kern\xgrayspace}}\kern\ygrayspace}}}
%
%
\long\def\surroundboxb#1{\leavevmode\hbox{\vtop{%
\vbox{\kern\ygrayspace%
\hbox{\kern\xgrayspace\vbox{\advance\hsize-2\xgrayspace#1}%
      \kern\xgrayspace}}\kern\ygrayspace}}}
%
%
%
\long\def\PScommands{%
\special{rawpostscript
/sharpbox{%
           newpath
           xmin ymin moveto
           xmin ymax lineto
           xmax ymax lineto
           xmax ymin lineto
           xmin ymin lineto
           closepath 
          }bind def
}%
\special{rawpostscript
/sharpboxnb{%
           newpath
           xmin ymin moveto
           xmin ymax lineto
           xmax ymax lineto
           xmax ymin lineto
          }bind def
}%
\special{rawpostscript
/sharpboxnt{%
           newpath
           xmin ymax moveto
           xmin ymin lineto
           xmax ymin lineto
           xmax ymax lineto
          }bind def
}%
\special{rawpostscript
/roundbox{%
           newpath
           xmin radius add ymin moveto
           xmax ymin xmax ymax radius arcto
           xmax ymax xmin ymax radius arcto
           xmin ymax xmin ymin radius arcto
           xmin ymin xmax ymin radius arcto 16 {pop} repeat
           closepath
          }bind def
}%
\special{rawpostscript
/sharpcorners{%
               sharpbox gsave grayshade setgray fill grestore 
               linewidth setlinewidth stroke
              }bind def
}%
\special{rawpostscript
/sharpcornersnt{%
               sharpboxnt gsave grayshade setgray fill grestore 
               linewidth setlinewidth stroke
              }bind def
}%
\special{rawpostscript
/sharpcornersnb{%
               sharpboxnb gsave grayshade setgray fill grestore 
               linewidth setlinewidth stroke
              }bind def
}%
\special{rawpostscript
/roundcorners{%
               roundbox gsave grayshade setgray fill grestore 
               linewidth setlinewidth stroke
              }bind def
}%
\special{rawpostscript
/plainbox{%
           sharpbox grayshade setgray fill 
          }bind def
}%
%
\special{rawpostscript
/roundnoframe{%
               roundbox grayshade setgray fill 
              }bind def
}%
\special{rawpostscript
/sharpnoframe{%
               sharpbox grayshade setgray fill 
              }bind def
}%
}%
%
%

\def\pshade#1{\Parashade{sharpcorners}{#1}{0.95}{10}{0.5}{10}{10}}


\def\boxit#1{\vbox{\hrule\hbox{\vrule\kern3pt
                                \vbox{\kern3pt#1\kern3pt}\kern3pt\vrule}\hrule}}

\def\boxitnb#1{\vbox{\hrule\hbox{\vrule\kern3pt
                                \vbox{\kern3pt#1\kern3pt}\kern3pt\vrule}}}

\def\boxitnt#1{\vbox{\hbox{\vrule\kern3pt
                                \vbox{\kern3pt#1\kern3pt}\kern3pt\vrule}\hrule}}

\long\def\boxtext#1#2{$$\boxit{\vbox{\hsize #1\noindent\strut #2\strut}}$$}



\def\texshopbox#1{\boxtext{462pt}{\vskip-1.5pc\pshade{\vskip-1.0pc#1\vskip-2.0pc}}}


%
%
%
%
%
%
%
%
\font\helbigbig=cmr10 scaled 2500%
\font\helbigb=cmbx10 scaled 1500%
\font\eightbold=cmbx8%

\def\tenf{\hel}%
\def\tenit{\heli}%
\def\ninef{\ninehel}%
\def\nineit{\nineheli}%
%
%


\font\tenrm=cmr10%
\font\teni=cmmi10%
\font\tensy=cmsy10%
\font\tenbf=cmbx10%
\font\tentt=cmtt10%
\font\tenit=cmti10%
\font\tensl=cmsl10%

\def\tenpoint{\def\rm{\fam0\tenrm}%
\textfont0=\tenrm%
\textfont1=\teni%
\textfont2=\tensy%
\textfont\itfam=\tenit%
\textfont\slfam=\tensl%
\textfont\ttfam=\tentt%
\textfont\bffam=\tenbf%
\scriptfont0=\sevenrm%
\scriptfont1=\seveni%
\scriptfont2=\sevensy%
\scriptscriptfont0=\sixrm%
\scriptscriptfont1=\sixi%
\scriptscriptfont2=\sixsy%
\def\it{\fam\itfam\tenit}%
\def\tt{\fam\ttfam\tentt}%
\def\sl{\fam\slfam\tensl}%
\scriptfont\bffam=\sevenbf%
\scriptscriptfont\bffam=\sixbf%
\def\bf{\fam\bffam\tenbf}%
\normalbaselineskip=18pt%
\normalbaselines\rm}%

\font\ninerm=cmr9%
\font\ninebf=cmbx9%
\font\nineit=cmti9%
\font\ninesy=cmsy9%
\font\ninei=cmmi9%
\font\ninett=cmtt9%
\font\ninesl=cmsl9%

\def\ninepoint{\def\rm{\fam0\ninerm}%
\textfont0=\ninerm%
\textfont1=\ninei%
\textfont2=\ninesy%
\textfont\itfam=\nineit%
\textfont\slfam=\ninesl%
\textfont\ttfam=\ninett%
\textfont\bffam=\ninebf%
\scriptfont0=\sixrm%
\scriptfont1=\sixi%
\scriptfont2=\sixsy%
\def\it{\fam\itfam\nineit}%
\def\tt{\fam\ttfam\ninett}%
\def\sl{\fam\slfam\ninesl}%
\scriptfont\bffam=\sixbf%
\scriptscriptfont\bffam=\fivebf%
\def\bf{\fam\bffam\ninebf}%
\normalbaselineskip=16pt%
\normalbaselines\rm}%

\font\eightrm=cmr8%
\font\eighti=cmmi8%
\font\eightsy=cmsy8%
\font\eightbf=cmbx8%
\font\eighttt=cmtt8%
\font\eightit=cmti8%
\font\eightsl=cmsl8%

\def\eightpoint{\def\rm{\fam0\eightrm}%
\textfont0=\eightrm%
\textfont1=\eighti%
\textfont2=\eightsy%
\textfont\itfam=\eightit%
\textfont\slfam=\eightsl%
\textfont\ttfam=\eighttt%
\textfont\bffam=\eightbf%
\scriptfont0=\sixrm%
\scriptfont1=\sixi%
\scriptfont2=\sixsy%
\scriptscriptfont0=\fiverm%
\scriptscriptfont1=\fivei%
\scriptscriptfont2=\fivesy%
\def\it{\fam\itfam\eightit}%
\def\tt{\fam\ttfam\eighttt}%
\def\sl{\fam\slfam\eightsl}%
\scriptscriptfont\bffam=\fivebf%
\def\bf{\fam\bffam\eightbf}%
\normalbaselineskip=14pt%
\normalbaselines\rm}%

\font\sevenrm=cmr7%
\font\seveni=cmmi7%
\font\sevensy=cmsy7%
\font\sevenbf=cmbx7%

\font\sixrm=cmr6%
\font\sixi=cmmi6%
\font\sixsy=cmsy6%
\font\sixbf=cmbx6%

\fontdimen13\tensy=2.6pt%
\fontdimen14\tensy=2.6pt%
\fontdimen15\tensy=2.6pt%
\fontdimen16\tensy=1.2pt%
\fontdimen17\tensy=1.2pt%
\fontdimen18\tensy=1.2pt%

\def\tenf{\tenpoint}%
\def\ninef{\ninepoint}%
%




\def\texshopbox#1{\boxtext{462pt}{\vskip-1.5pc\pshade{\vskip-1.0pc#1\vskip-2.0pc}}}


\input miniltx

\ifx\pdfoutput\undefined
  \def\Gin@driver{dvips.def} 
\else
  \def\Gin@driver{pdftex.def} 
\fi

\input graphicx.sty
\resetatcatcode

\long\def\fig#1#2#3{\vbox{\vskip1pc\vskip#1
\prevdepth=12pt \baselineskip=12pt
\vskip1pc
\hbox to\hsize{\hfill\vtop{\hsize=30pc\noindent{\eightbf Figure #2\ }
{\eightpoint#3}}\hfill}}}

\def\show#1{}

\def\frac#1#2{{#1\over #2}}

\rightheadline{\botmark}

\pageno=1

\rightheadline{\botmark}

\pn {\bf February 2020}
\bigskip \bigskip \bigskip \bigskip

\bigskip\bigskip\bigskip

\def\longpapertitle#1#2#3{{\bf \centerline{\helbigb
{#1}}}\medskip{\bf \centerline{\helbigb
{#2}}}\medskip{\centerline{
by}}\medskip{\bf \centerline{
{#3}}}\bigskip}

\longpapertitle{Constrained Multiagent Rollout and}{Multidimensional Assignment with the Auction Algorithm}{{Dimitri Bertsekas\footnote{\dag}{\ninepoint McAfee Professor of Engineering, MIT, Cambridge, MA, and Fulton Professor of Computational Decision Making, ASU, Tempe, AZ.}}}

\centerline{\bf Abstract}

We consider an extension of
the rollout algorithm that applies to constrained deterministic dynamic programming, including challenging combinatorial optimization problems. The algorithm relies on a suboptimal policy, called base heuristic. Under suitable assumptions, we show that if the base heuristic
produces a feasible solution, the rollout algorithm has a cost improvement property: it produces a feasible solution, whose cost is no worse
than the base heuristic's cost. 

We then focus on multiagent problems, where the control at each stage consists of multiple components (one per agent), which are coupled either through the cost function or the constraints or both. We show that the cost improvement property is maintained with an alternative implementation that has greatly reduced computational requirements, and makes possible the use of rollout in problems with many agents. We demonstrate this alternative algorithm by applying it to layered graph problems that involve both a spatial and a temporal structure. We consider in some detail a prominent example of such problems: multidimensional assignment, where we use the auction algorithm for 2-dimensional assignment as a base heuristic. This auction algorithm is particularly well-suited for our context, because through the use of prices, it can advantageously use the solution of an assignment problem as a starting point for solving other related assignment problems, and this can greatly speed up the execution of the rollout algorithm.

\vfill\eject

\section{Introduction}
\mark{Introduction}

\vskip-1pc

\pn We consider a deterministic optimal control problem
involving the system
$$x_{k+1}=f_k(x_k,u_k),\qquad k=0,\ldots,N-1,$$ 
where $x_k$ is the state, taking values in some (possibly infinite) set, $u_k$ is the control at time $k$, taking values in some finite set, and $f_k$ is some function. The initial state is given and is denoted by $x_0$. A
sequence of the form
$$T=(x_0,u_0,x_1,u_1,\ldots,u_{N-1},x_N),\xdef\systraject{\lab}\eqnum\show{spconst}$$
where 
$$x_{k+1}=f_k(x_k,u_k),\qquad k=0,1,\ldots,N-1,\xdef\sysequation{\lab}\eqnum\show{spconst}$$
is referred to as a {\it complete trajectory\/}. We distinguish a complete trajectory from a {\it partial
trajectory\/}, which is defined to be a subset of a complete trajectory, consisting of a subsequence of time-contiguous states and controls. 
Our problem is stated succinctly as
$$\min_{T\in C} G(T),\xdef\multicostone{\lab}\eqnum\show{spconst}$$
where $G$ is a given real-valued cost function and $C$ is a given constraint set of trajectories.\footnote{\dag}{\ninepoint Actually it is not essential that we know the form of the function $G$. Instead it is sufficient to have access to a human or software expert that enables us to compare any two trajectories $T_1$ and $T_2$, without assigning numerical values to them. It is essential, however, that the expert's rankings should have a transitivity property: if $T$ is ranked better than $T'$ and $T'$ is ranked better than $T''$, then $T$ is ranked better than $T''$. Of course, the expert should also be able to determine whether a given trajectory $T$ satisfies the constraint $T\in C$.}

As an example, we note the common special case of the additive cost
$$G(x_0,u_0,x_1,u_1,\ldots,u_{N-1},x_N)=g_N(x_N)+\sum_{k=0}^{N-1}g_k(x_k,u_k),\xdef\addcost{\lab}\eqnum\show{spconst}$$ 
where $g_k$, $k=0,1,\ldots, N$, are given real-valued functions, and the controls satisfy the time-uncoupled (but state-dependent) constraints 
$$u_k\in U_k(x_k),\qquad k=0,1,\ldots, N-1,\xdef\uncoupledcontrol{\lab}\eqnum\show{spconst}$$
(so here $C$ is the set of trajectories that are generated by the system equation with controls satisfying the above constraints).
This is a standard problem formulation, which is usually taken as the starting point for deterministic dynamic programming (DP for short). Our aim, however, is to address problems involving far more complicated constraints, for which the exact solution of the problem is typically intractable, including multiagent problems for which $u_k$ consists of multiple components (one per agent). To this end, we will consider approximate solution methods based on the rollout approach.

The general idea of a rollout algorithm is to start with a suboptimal solution method called {\it base heuristic\/}, and to aim at {\it cost improvement\/}: a guarantee (under suitable assumptions) that the rollout algorithm produces a feasible solution, whose cost is no worse than the cost corresponding to the base heuristic. Rollout is simple and reliable, and has been used with considerable success for general unconstrained deterministic DP
problems of the form \multicostone-\uncoupledcontrol. It has also been used in  minimax/game settings, and in stochastic settings, in conjunction with simulation-based evaluation of the expected cost produced by the base heuristic starting from a given state. For applications and related work we refer to Tesauro and Galperin [TeG96], Bertsekas and Tsitsiklis [Ber96], Bertsekas, Tsitsiklis, and Wu [BTW97], 
Bertsekas [Ber97], Christodouleas [Chr97], Bertsekas and  Casta\~ non [BeC99], Duin and Voss [DuV99],
Secomandi [Sec00], [Sec01], [Sec03], Ferris and Voelker [FeV02], [FeV04],  McGovern, Moss, and Barto [MMB02], Savagaonkar, Givan, and
Chong [SGC02], Guerriero and Mancini [GuM03], Tu and Pattipati [TuP03], Wu, Chong, and Givan [WCG03], Chang, Givan, and Chong [CGC04],  Meloni, Pacciarelli, and Pranzo [MPP04], Yan, Diaconis, Rusmevichientong, and Van Roy [YDR04], Bertsekas [Ber05a], [Ber05b], Besse and Chaib-draa [BeC08], Sun et al.\ [SZL08], Bertazzi et al.\ [BBG13], Sun et al.\ [SLJ13], Tesauro et al.\ [TGL13],  Antunes and Heemels [AnH14], Beyme and Leung [BeL14], Goodson, Thomas, and Ohlmann [GTO15], Khashooei, Antunes, and Heemels [KAH15], Li and Womer [LiW15], Mastin and Jaillet [MaJ15], Huang,  Jia, and Guan [HJG16], Simroth, Holfeld, and Brunsch [SHB15], Lan, Guan, and Wu [LGW16], Ulmer [Ulm17], Bertazzi and Secomandi [BeS18], Guerriero, Di Puglia, and Macrina [GDM19], Sarkale et al.\ [SNC18], Ulmer at al.\ [UGM18], Bertsekas [Ber19c], and Chu, Xu, and Li [CXL19]. These works discuss variants and problem-specific adaptations of rollout algorithms for a broad variety of practical problems, and consistently report  positive computational experience.

In this paper we will adapt the rollout approach to construct methods that can address suboptimally the constrained DP problem \multicostone. Our line of analysis and development are based on the ideas of the paper by Bertsekas, Tsitsiklis, and Wu [BTW97], where rollout was applied to general discrete deterministic optimization problems. Related constrained rollout ideas are also discussed in the author's papers [Ber05a], [Ber05b]. The extension to constrained multiagent rollout (see Section 4) is based on the  author's recent paper [Ber19b], which proposed a modification of the standard rollout algorithm to deal efficiently with the special computational demands of the many-agent case.

Generally constrained DP problems can be transformed to unconstrained DP problems. The idea is to redefine the state at stage $k$ to be the partial trajectory
$$y_k=(x_0,u_0,x_1,\ldots, u_{k-1}, x_k),$$
which evolves according to a redefined system equation:
$$y_{k+1}=\br(y_k,u_k,f_k(x_k,u_k)\br).\xdef\trajsystem{\lab}\eqnum\show{spconst}$$
The problem then becomes to find a control sequence that minimizes the terminal cost
$G(y_N)$ subject to the constraint $y_N\in C$. This is a problem to which the standard form of DP applies. 

Unfortunately, with the DP reformulation just described, the exact solution of the problem is typically impractical because the
associated computation can be overwhelming. It is much greater than the computation for the
corresponding additive cost/time-uncoupled control constraints problem \addcost-\uncoupledcontrol, where the constraint
$T\in C$ is absent. This is true even when $C$ is specified in terms
of a finite number of constraint functions that are time-additive, i.e., $T\in C$ if 
$$g_N^m(x_N)+\sum_{k=0}^{N-1}g_k^m(x_k,u_k)\le
b^m,\qquad m=1,\ldots,M,\xdef\multiconstradd{\lab}\eqnum\show{spconst}$$
where  $g_k^m$, $k=0,1,\ldots, N$, and $b^m$, $m=1,\ldots,M$, are given functions and scalars, respectively.
For further appreciation of the issues involved, the reader may consult the author's textbook [Ber17], and the journal literature, which contains several proposals for suboptimal
solution of the problem in the case where the constraints are of the form \multiconstradd, using among others, multiobjective optimization ideas; see e.g., Jaffe
[Jaf84],  Martins [Mar84], Guerriero and Musmanno [GuM01],  and Stewart and White [StW91], who also survey
earlier work. Generally, experience with constrained DP problems suggests that the use of an approximate solution approach is essentially unavoidable. This is the motivation for the methodology of this paper.

\subsection{Using a Base Heuristic for  Constrained Rollout}

\pn We will now describe our rollout algorithm. We assume the availability of a base heuristic, which for any given partial trajectory
$$y_k=(x_0,u_0,x_1,\ldots, u_{k-1}, x_k),$$
can produce a (complementary) partial  trajectory\footnote{\dag}{\ninepoint The nature of the base heuristic is essentially arbitrary, and may strongly depend on $y_k$ as well as $k$. For an extreme but practically interesting possibility, we may have a partition of the set of partial trajectories $y_k$, and a collection of multiple heuristics that are specially tailored to the sets of the partition. We may then select the appropriate heuristic to use on each set of the partition. Of course the properties of the rollout algorithm and its potential for cost improvement will depend on the design of the base heuristic, and in Section 2, we will discuss properties that are favorable in this regard.} 
$$R(y_k)=(x_k, u_k,x_{k+1},u_{k+1},\ldots,u_{N-1},x_N),$$
that starts at $x_k$ and satisfies for every $t=k,\ldots,N-1$ the system equation
$$x_{t+1}=f_t(x_t,u_t).$$
Thus, given $y_k$ and any control $u_k$, we can use the base heuristic to obtain a complete trajectory as follows:
\nitem{(a)} Generate the next state $x_{k+1}=f_k(x_k,u_k)$.
\nitem{(b)} Extend $y_k$ to obtain the partial trajectory
$$y_{k+1}=\big(y_k,u_k,f_k(x_k,u_k)\big).$$
\nitem{(c)} Run the base heuristic from $y_{k+1}$ to obtain the partial trajectory $R(y_{k+1})$. 
\nitem{(d)} Join the two partial trajectories $y_{k+1}$ and $R(y_{k+1})$ to obtain the  complete trajectory $\big(y_k,u_k,R(y_{k+1})\big)$, which is denoted by $T_k(y_k,u_k)$:
$$T_k(y_k,u_k)=\big(y_k,u_k,R(y_{k+1})\big).\xdef\completetr{\lab}\eqnum\show{spconst}$$
\smskip

A complete trajectory $T_k(y_k,u_k)$ of the form \completetr\ is generally feasible for only a subset of controls $u_k$, which we denote by $U_k(y_k)$:
$$U_k(y_k)=\big\{u_k\mid T_k(y_k,u_k)\in C\big\}.\xdef\futconstr{\lab}\eqnum\show{spconst}$$
Our rollout algorithm starts from a given initial state $\tl y_0=\tl x_0$, and generates successive partial trajectories $\tl y_1,\ldots,\tl y_N$, of the form
$$\tl y_{k+1}=\big(\tl y_k,\tl u_k,f_k(\tl x_k,\tl u_k)\big),\qquad k=0,\ldots,N-1,\xdef\systemevolve{\lab}\eqnum\show{spconst}$$
where $\tl x_k$ is the last state component of $\tl y_k$, and $\tl u_k$ is a control that minimizes the heuristic cost $H\big(T_k(\tl y_k,u_k)\big)$ over all  $u_k$ for which $T_k(\tl y_k,u_k)$ is feasible. Thus at stage $k$,
the algorithm forms  the set $U_k(\tl y_k)$ and selects from $U_k(\tl y_k)$ a
control
$\tl u_k$ that minimizes  the cost of the complete trajectory $T_k(\tl y_k,u_k)$:
$$\tl u_k\in\arg\min_{u_k\in U_k(\tl y_k)}G\big(T_k(\tl y_k,u_k)\big);\xdef\mincontr{\lab}\eqnum\show{spconst}$$ 
see Fig.\ 1.1.
The objective is to produce a feasible final complete trajectory $\tl y_N$, which has a cost $G(\tl y_N)$ that is no larger than the cost of $R(\tl y_0)$ produced by the base heuristic starting from $\tl y_0$, i.e., 
$$G(\tl y_N)\le G\big(R(\tl y_0)\big).\xdef\costimprove{\lab}\eqnum\show{spconst}$$

\topinsert
\centerline{\hskip0pc\includegraphics[width=4.8in]{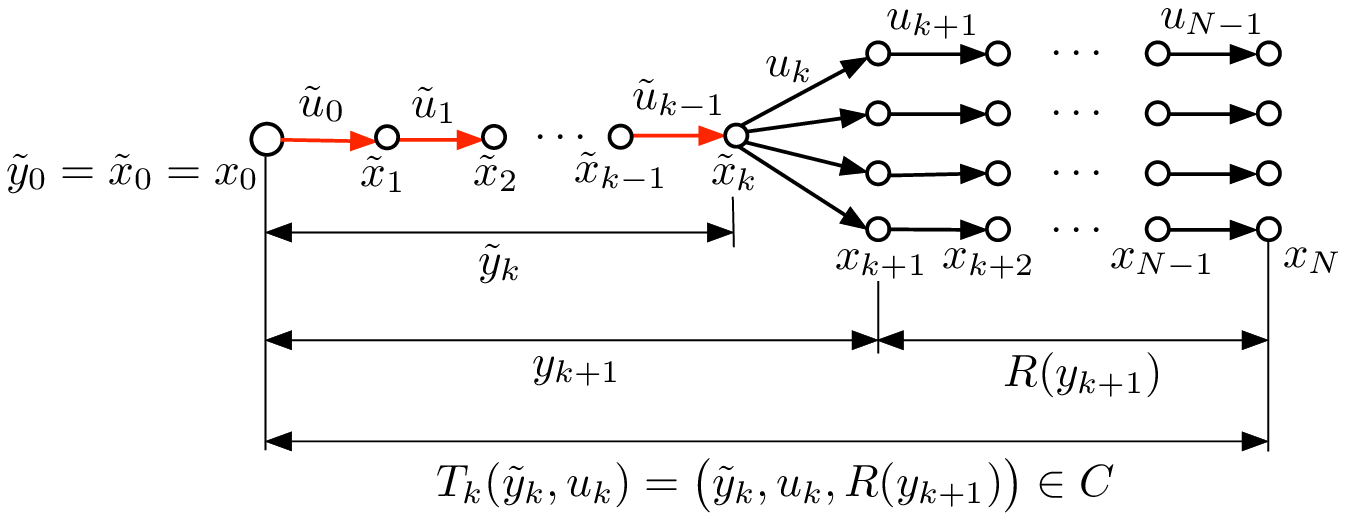}}
\vskip-1pc
\hskip-2pc\fig{0pc}{1.1.} {Illustration of the rollout algorithm. At stage $k$, and given the current partial trajectory
$$\tl y_k=(\tl x_0,\tl u_0,\tl x_1,\ldots,\tl u_{k-1},\tl x_k),$$ 
which starts at $\tl x_0$ and ends at $\tl x_k$, it considers all possible next states 
$x_{k+1}=f_k(\tl x_k,u_k)$, and runs the base heuristic starting at $y_{k+1}=(\tl y_k,u_k,x_{k+1})$. It then:
\nitem{(a)} Finds $\tl u_k$, the
control that minimizes over all $u_k$ the cost $G\big(T_k(\tl y_k,u_k)\big)$ over all $u_k$ that yield a feasible complete trajectory $T_k(\tl y_k,u_k)\in C$.
\nitem{(b)} Extends $\tl y_k$ by $\big(\tl u_k, f_k(\tl x_k,\tl u_k)\big)$ to form $\tl y_{k+1}$.
}\endinsert

Note that $T_k(\tl y_k,u_k)$ is not guaranteed to be feasible for any given $u_k$ (i.e., it may not belong to $C$), but for the analysis of Section 2, our assumptions will guarantee that the constraint set $U_k(\tl y_k)$
of the problem \mincontr\ is nonempty, so that our rollout algorithm is well-defined. Later, in Section 3, we will modify our algorithm so that it is well-defined under the weaker assumption that just {\it the complete trajectory generated by the base heuristic starting from the given initial state $\tl y_0$ is feasible\/}, i.e., $R(\tl y_0)\in C.$

It can be seen that our constrained rollout algorithm is not much more complicated or
computationally demanding than its unconstrained version where the constraint $T\in C$ is not present (as long as checking feasibility of
a complete trajectory $T$ is not computationally demanding). Note, however, that our algorithm makes essential use of the deterministic character of the problem,
and does not admit a straightforward extension to stochastic or minimax problems, since checking feasibility of a complete trajectory is typically difficult in stochastic and minimax contexts. 

In the next section, we establish the cost improvement
properties of our algorithm under the assumptions of sequential consistency and sequential improvement introduced in [BTW97]. In Section 3, we discuss some variants and extensions of the algorithm. In Section 4, we discuss how the algorithm can be applied to multiagent problems with computational complexity that is linear in the number of agents, rather than exponential. In Section 5, we discuss the application of our algorithm to combinatorial problems, including a broad class of layered graph problems, which contains as special cases problems of scheduling, routing, and multidimensional assignment. We will illustrate rollout primarily through the multidimensional assignment problem by using base heuristics that solve repeatedly 2-dimensional assignment problems with the auction algorithm. This algorithm is particularly well-suited to our context because it allows the efficient use of the solution of a given problem as a starting point for the solution of a related assignment problem. Detailed presentations of the auction algorithm are given in several sources, including
the author's textbooks [Ber91], [Ber98], and the survey paper [Ber92]. For purposes of easy reference and completeness, we provide an appendix with a description, which is based on the author's brief tutorial paper [Ber01].

\vskip-1pc

\section{Cost Improvement with the Rollout Algorithm}
\mark{Cost Improvement with the Rollout Algorithm}

\vskip-1pc

\pn Let us summarize the rollout algorithm \systemevolve-\mincontr. It starts at stage 0 and
sequentially proceeds to the last stage. At stage $k$, it maintains  a partial trajectory 
$$\tl y_k=(\tl x_0,\tl u_0,\tl x_1,\ldots,\tl u_{k-1},\tl x_k)\xdef\partialtraject{\lab}\eqnum\show{spconst}$$ 
that
starts at the given initial state
$\tl y_0=\tl x_0$, and is such that 
$$\tl x_{t+1}=f_t(\tl x_t,\tl u_t),\qquad t=0,1,\ldots,k-1.$$
The algorithm then forms the set of controls 
$$U_k(\tl y_k)=\big\{u_k\mid T_k(\tl y_k,u_k)\in C\big\}$$
that is consistent with feasibility [cf.\ Eq.\ \futconstr], and chooses a control $\tl u_k\in U_k(\tl y_k)$ according to the minimization 
$$\tl u_k\in\arg\min_{u_k\in U_k(\tl y_k)}G\big(T_k(\tl y_k,u_k)\big),\xdef\mincontrt{\lab}\eqnum\show{spconst}$$ 
[cf.\ Eq.\ \mincontr],
 where
 $$T_k(\tl y_k,u_k)=\Big(\tl y_k,u_k, R\big(\tl y_k,u_k, f_k(\tl x_k,u_k)\big)\Big);$$
 [cf.\ Eq.\ \completetr]. Finally, the algorithm sets  
 $$\tl x_{k+1}=f_k(\tl x_k,\tl u_k),\qquad \tl y_{k+1}=(\tl y_k,\tl u_k,\tl x_{k+1}),\eqnum\show{spconst}$$
[cf.\ Eq.\ \systemevolve]. 

We will introduce conditions guaranteeing that the control set $U_k(\tl y_k)$ in the minimization \mincontrt\ is nonempty, and that the costs of the complete trajectories $T_k(\tl y_k,\tl u_k)$ are improving with each $k$ in the sense that
$$G\big(T_{k+1}(\tl y_{k+1},\tl u_{k+1})\big)\le G\big(T_{k}(\tl y_{k},\tl u_{k})\big),\qquad k=0,1,\ldots,N-1,$$
while at the first step of the algorithm we have
$$G\big(T_{0}(\tl y_{0},\tl u_{0})\big)\le G\big(R(\tl y_0)\big).$$
It will then follow that the cost improvement condition $G(\tl y_N)\le G\big(R(\tl y_0)\big)$ [cf.\ Eq.\ \costimprove] holds.

\texshopbox{\definition{2.1:} We say that the base heuristic is {\it sequentially consistent} if whenever it
generates a partial trajectory 
$$(x_k,u_k,x_{k+1},u_{k+1},\ldots,u_{N-1},x_N),$$
starting from a partial trajectory $y_k$, it also generates the  partial trajectory 
$$(x_{k+1},u_{k+1},x_{k+2},u_{k+2},\ldots,u_{N-1},x_N),$$
starting from the partial trajectory $y_{k+1}=\big(y_k,u_k,x_{k+1}\big)$.
} 

Sequentially consistent heuristics are often used in practice. For example greedy heuristics tend to be sequentially consistent. Also any policy [a sequence of feedback control functions $\m_k(y_k)$, $k=0,1,\ldots,N-1$] for the DP problem of minimizing the terminal cost $G(y_N)$ subject to the system equation $y_{k+1}=\big(y_k,u_k,f_k(x_k,u_k)\big)$ and the feasibility constraint $y_N\in C$ [cf.\ Eq.\ \trajsystem] can be seen to be sequentially consistent. 

For a given partial trajectory $y_k$, let us denote by $y_k\cup R(y_k)$ the complete trajectory obtained by joining $y_k$ with the partial trajectory generated by the base heuristic starting from $y_k$. Thus if  $y_k=(x_0,u_0,\ldots,u_{k-1},x_k)$ and $R(y_k)=(x_k,u_{k+1},\ldots,u_{N-1},x_N)$, we have
$$y_k\cup R(y_k)=(x_0,u_0,\ldots,u_{k-1},x_k,u_{k+1},\ldots,u_{N-1},x_N).$$

\texshopbox{
\definition{2.2:} We say that the base heuristic is {\it sequentially improving} if for every $k$ and partial trajectory $y_k$ for which $y_k\cup R(y_k)\in C$, the set $U_k(y_k)$ is nonempty, and we have 
$$ G\big(y_k\cup R(y_k)\big)\ge \min_{u_k\in U_k(y_k)}G\big(T_k(y_k,u_k)\big).\xdef\seqimpo{\lab}\eqnum\show{spconst}$$
} 

Note that if the base heuristic is sequentially consistent, it is
also sequentially improving. The reason is that for a sequentially consistent heuristic, $y_k\cup R(y_k)$ is equal to one of the trajectories contained in the set $\big\{T_k(y_k,u_k)\mid u_k\in U_k(y_k)\big\}$.

Our main result is contained in the following proposition. \old{To state compactly the
result, consider the  partial
 trajectory 
$$\tl y_k=(\tl x_0,\tl u_0,\tl x_1,\ldots,\tl u_{k-1},\tl x_k)$$
maintained by the rollout
algorithm after $k$ stages, the set of trajectories 
$$\bl\{T_k(\tl y_k,u_k)\mid u_k\in U_k(\tl y_k)\br\}$$ 
that are feasible, and the complete trajectory $T_k(\tl y_k,\tl u_k)$ within this set that corresponds to the control
$\tl u_k$ chosen by the rollout algorithm, i.e.,
$$T_k(\tl y_k,\tl u_k)=\big(\tl y_k,\tl u_k, R(\tl y_{k+1})\big),\eqnum\show{spconst}$$
where
$$\tl y_{k+1}=\big(\tl y_k,\tl u_k,f_k(\tl x_k,\tl u_k)\big).$$}

\texshopbox{
\proposition{2.1:} Assume that the base heuristic is sequentially improving and generates a feasible complete trajectory starting from the initial state $\tl y_0=\tl x_0$, i.e., $R(\tl y_0)\in C$. Then for each $k$,
the set $U_k(\tl y_k)$ is nonempty, and we have
$$G\bl(R(\tl y_0)\br)\ge G\big(T_0(\tl y_0,\tl u_0)\big)\ge G\big(T_1(\tl y_1,\tl u_1)\big)\ge  \cdots\ge G\big(T_{N-1}(\tl y_{N-1},\tl u_{N-1})\big)=  G(\tl y_N),$$
where 
$$T_k(\tl y_k,\tl u_k)=\big(\tl y_k,\tl u_k,R(\tl y_{k+1})\big);$$
cf.\ Eq.\ \completetr.  In particular, the final trajectory $\tl y_N$ generated by the constrained rollout algorithm is feasible and has no larger cost than the trajectory $R(\tl y_0)$ generated by the base heuristic starting from the initial state.
}

\proof Consider $R(\tl y_0)$, the complete trajectory generated by the base heuristic starting from $\tl y_0$. Since $\tl y_0\cup R(\tl y_0)=R(\tl y_0)\in C$ by assumption, it follows from the sequential improvement definition, that the set $U_0(\tl y_0)$ is nonempty and we have
$$G\bl(R(\tl y_0)\br)\ge G\big(T_0(\tl y_0,\tl u_0)\big),$$
[cf.\ Eq.\ \seqimpo], while $T_0(\tl y_0,\tl u_0)\in C$.

The preceding argument can be repeated for the next stage, by replacing $\tl y_0$ with $\tl y_1$, and $R(\tl y_0)$ with
$T_0(\tl y_0,\tl u_0)$. Since $\tl y_1\cup R(\tl y_1)=T_0(\tl y_0,\tl u_0)\in C$, from the sequential improvement definition, the set $U_1(\tl y_1)$ is nonempty and we have
$$G\bl(T_0(\tl y_0,\tl u_0)\br)= G\bl(\tl y_1\cup R(\tl y_1)\br)\ge G\bl(T_1(\tl y_1,\tl u_1)\br),$$ 
[cf.\ Eq.\ \seqimpo], while $T_1(\tl y_1,\tl u_1)\in C$.
Similarly, the argument can be successively repeated for every $k$, to verify  that $U_k(\tl y_k)$ is nonempty and that $G\bl(T_k(\tl y_k,\tl u_k)\br)\ge G\bl(T_{k+1}(\tl y_{k+1},\tl u_{k+1})\br)$ for all $k$. \qed

Proposition 2.1 implies that for a base heuristic that is sequentially improving and produces a feasible initial complete trajectory, starting from the initial state $\tl y_0$, the rollout algorithm generates at each stage $k$ a feasible complete 
trajectory that is no worse than its predecessor in terms of cost. It follows that the algorithm 
produces a  final complete trajectory $\tl y_N=(\tl x_0,\tl u_0,\tl x_1,\ldots,\tl u_{N-1},\tl x_N)$ that is feasible and has cost that
is no larger than the cost of the initial complete trajectory produced by the base heuristic. On the other hand it is easy to construct examples where the sequential improvement condition \seqimpo\ is violated and the cost of the solution produced by rollout is larger than the cost of the solution produced by the base heuristic starting from the initial state (see [Ber19a], Example 2.4.2).

\vskip-1pc

\section{Variants and Extensions of the Rollout Algorithm}
\mark{Variants and Extensions of the Rollout Algorithm}

\vskip-1pc

\pn We will now discuss some variations and extensions of the rollout algorithm of Section 2.

\subsection{Rollout Without the Sequential Improvement Assumption - Fortified Rollout}

\pn Let us consider the case where
the sequential improvement assumption is not satisfied. Then it may happen that given the current partial  trajectory $\tl y_k$, the set of controls $U_k(\tl y_k)$ 
that corresponds to feasible trajectories $T_k(\tl y_k,u_k)$ [cf.\ Eq.\ \futconstr] is empty, in which case the
rollout algorithm cannot extend the partial trajectory
$\tl y_k$ further. To bypass this difficulty, we propose a modification, called {\it fortified
rollout algorithm\/}, which is patterned after an algorithm given in [BTW97] for the
case of an unconstrained DP problem (see also [Ber17], Section 6.4, and [Ber19a], Section 2.4.1). For validity of this algorithm, {\it we require that the base heuristic generates a feasible complete 
trajectory $R(\tl y_0)$ starting from the initial state
$\tl y_0$\/}.

The fortified rollout algorithm, in addition to the current partial trajectory
$$\tl y_k=(\tl x_0,\tl u_0,\tl x_1,\ldots,\tl u_{k-1},\tl x_k),$$ 
maintains a complete trajectory
$\hat T_k$, called {\it tentative best trajectory\/}, which is feasible (i.e., $\hat T_k\in C$) and agrees with $\tl y_k$ up to state $\tl x_k$, i.e., $\hat T_k$ has the form
$$\hat T_k=(\tl x_0,\tl u_0,\tl x_1,\ldots,\tl u_{k-1},\tl x_k,\ol u_k,\ol x_{k+1},\ldots,\ol u_{N-1},\ol x_N),\xdef\tentative{\lab}\eqnum\show{spconst}$$
for some $\ol u_k,\ol x_{k+1},\ldots,\ol u_{N-1},\ol x_N$ such that 
$$\ol x_{k+1}=f_k(\tl x_k,\ol u_k),\qquad\ \ \ol x_{t+1}=f_t(\ol x_t,\ol u_t),\quad t=k+1,\ldots,N-1.$$ Initially,
$\hat T_0$ is the complete trajectory $R(\tl y_0)$, generated by the base heuristic starting from $\tl y_0$, which is assumed to be feasible. At stage
$k$, the algorithm  forms  the subset $\hat U_k(\tl y_k)$ of controls $u_k\in U_k(\tl y_k)$ such that the corresponding $T_k(\tl y_k,u_k)$ is not only feasible, but also has cost that is  no larger than the one of the current tentative best trajectory: 
$$\hat U_k(\tl y_k)=\Big\{u_k\in U_k(\tl y_k)\mid G\big(T_k(\tl y_k,u_k)\big)\le G(\hat T_k)\Big\}.$$
 
There are two cases to consider at state $k$:

\nitem{(1)} {\it The set $\hat U_k(\tl y_k)$ is nonempty\/}. Then the algorithm forms the partial trajectory  
$\tl y_{k+1}=(\tl y_k,\tl u_k,\tl x_{k+1}),$
where 
$$\tl u_k\in\arg\min_{u_k\in \hat U_k(\tl y_k)}G\big(T_k(\tl y_k,u_k)\big),\qquad \tl x_{k+1}=f_k(\tl x_k,\tl u_k),$$
and sets $T_k(\tl y_k,\tl u_k)$ as the new tentative best trajectory, i.e.,
$$\hat T_{k+1}=T_k(\tl y_k,\tl u_k).$$ 

\nitem{(2)} {\it The set $\hat U_k(\tl y_k)$ is empty\/}. Then, the algorithm forms the partial trajectory 
$\tl y_{k+1}=\big(\tl y_k,\tl u_k,\tl x_{k+1}),$ 
where 
$$\tl u_k=\ol u_k,\qquad \tl x_{k+1}=\ol x_{k+1},$$
and $\ol u_k,\ol x_{k+1}$ are the control and state subsequent to $\tl x_k$ in the current tentative best trajectory $\hat T_k$ [cf.\ Eq.\ \tentative], 
and leaves $\hat T_k$ unchanged, i.e., 
$$ \hat T_{k+1}=\hat T_k.$$

\smskip

It can be seen that the fortified rollout algorithm will follow the initial complete trajectory
$\tl T_0$, the one generated by the base heuristic starting from $\tl y_0$, up to a stage $\ol k$ where it will
discover a new feasible complete  trajectory with smaller cost to replace $\tl T_0$ as the tentative best trajectory.  Similarly, the new tentative best trajectory $\tl T_{\ol k}$
may be subsequently replaced by another feasible trajectory with smaller cost, etc. Note that if the
base heuristic is sequentially improving, and the fortified rollout algorithm will generate the same complete trajectory as the
(nonfortified) rollout algorithm given earlier, with the tentative best trajectory $\hat T_{k+1}$ being equal to the complete trajectory $T_k(\tl y_k,\tl u_k)$ for all $k$. The reason is that if the
base heuristic is sequentially improving the controls $\tl u_k$ generated by the nonfortified algorithm belong to the set $\hat U_k(\tl y_k)$ [by Prop.\ 2.1, case (1) above will hold]. 

However, it can be verified that even when the base heuristic is not sequentially improving, the fortified rollout
algorithm will generate a complete  trajectory that is feasible and has cost that is no worse than the cost of the
complete trajectory generated by the base heuristic starting from $\tl y_0$. This is because each tentative best trajectory has a cost that is no worse than the one of its predecessor, and the initial tentative best trajectory is just the trajectory generated by the base heuristic starting from the initial condition $\tl y_0$. 

\subsection{Tree-Based Rollout Algorithm}

\pn It is possible to improve the performance of the rollout algorithm at the expense of maintaining more than one
partial trajectory. In particular, instead of the partial trajectory $\tl y_k$ of Eq.\ \partialtraject, we can maintain a
{\it tree} of partial trajectories that is rooted at $\tl y_0$. These trajectories need not be of equal length, i.e., they
need not have the same number of stages. At each step of the algorithm, we select a single partial trajectory from this
tree, and execute the rollout algorithm's step as if this partial trajectory were the only one. Let this partial trajectory have $k$ stages and denote it by
$\tl y_k$. Then we extend
$\tl y_k$ similar to the rollout algorithm of Section 2, with possibly multiple feasible trajectories. 
There is also a fortified version of this algorithm where a tentative best trajectory is maintained, which is the minimum cost complete trajectory generated thus far. 

The aim of the tree-based algorithm is to obtain
improved performance, essentially because it can return to extend partial trajectories that were generated and temporarily abandoned at previous stages. The net result is a more flexible algorithm that is capable of examining more alternative trajectories. Note also that there is considerable
freedom to select the number of partial trajectories maintained in the tree. 

We finally mention a drawback of the tree-based algorithm: it is suitable for off-line computation, but it cannot be
applied in an on-line context, where the rollout control selection is made after the
current state becomes known as the system evolves in real-time. By contrast, the rollout algorithm of Section 2 and its fortified version are well-suited for on-line application.

\subsection{Rollout Without States}

\pn It is important to note that the constrained deterministic optimal control problem of this paper is very general. In particular, it contains as a special case the fully unstructured discrete optimization problem:
$$\eqalign{&\hbox{minimize\ \ }G(u)\cr &\hbox{subject to\ \
}u\in C,\cr}\xdef\gendiscrete{\lab}\eqnum\show{spconst}$$
where each solution $u$ has $N$ components; i.e., it has the form $u=(u_0,\ldots,u_{N-1})$, where $N$ is a positive integer, $C$ is a finite set of feasible solutions,  and $G(u)$ is some cost function.\footnote{\dag}{\ninepoint The reverse is also true, namely that any constrained deterministic optimal control problem of the form  \systraject-\multicostone, can be converted to the general discrete optimization form \gendiscrete, simply by expressing the states $x_k$ as functions of the preceding controls $u_0,\ldots,u_{k-1}$ through the system equation \sysequation, and eliminating them from the cost function expression and the constraints. This abstraction of the problem may be of value in some contexts because of its inherent simplicity.}  This is simply the special case of the deterministic optimal control problem where {\it each state $x_k$ can only take a single value}. Then the state space for each $k$ has a single element, and the system equation $x_{k+1}=f_k(x_k,u_k)$ is trivial and superfluous. Then in effect the partial trajectory $y_k$ is the
$k$-tuple $(u_0,\ldots,u_{k-1})$  consisting of the first $k$ components of a solution. 

We associate such a $k$-tuple with the $k$th stage of the finite horizon DP problem shown in Fig.\ 3.1. In
particular,  for $k=0,\ldots,N-1$, we view as the states of the $k$th stage all the possible $k$-tuples
$(u_0,\ldots,u_{k-1})$. The initial state is some artificial state. From this state
we may move to any state
$(u_0)$, with
$u_{0}$ belonging to the  set 
$$U_{0}=\bl\{u_0\mid \hbox{there exists a solution of the form }(
u_0,\ol u_1,\ldots,\ol u_{N-1})\in U\br\}.$$ 
Thus $U_0$ is the set of choices of $u_0$ that are consistent with feasibility.

\topinsert
\centerline{\hskip0pc\includegraphics[width=5.0in]{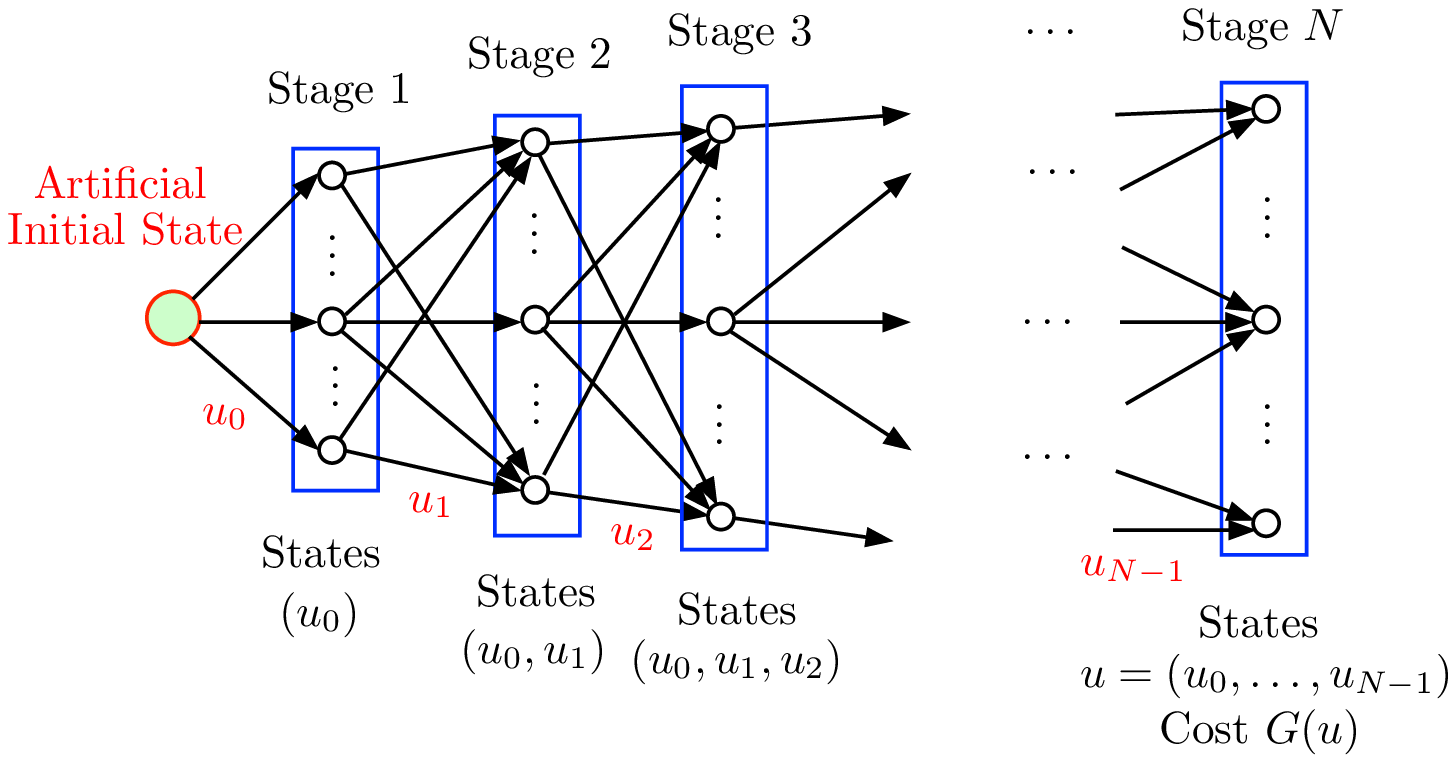}}
\vskip-1pc
\hskip-2pc\fig{0pc}{3.1.} {Formulation of a discrete optimization problem as a DP problem with $N$ stages. Starting with an artificial initial state, a new control component is chosen at each stage. There is a cost $G(u)$ only at the terminal stage.}\endinsert

More generally, from
a state
$(u_0,\ldots,u_{k-1}),$
we may move to any state of the form
$(u_0,\ldots,u_{k-1},u_{k}),$
such that $u_{k}$ belongs to the set
$$U_{k}(u_0,\ldots,u_{k-1})=\big\{u_{k}\mid\ \hbox{there exists a solution of the form } (u_0,\ldots,u_{k},\ol u_{k+1},\ldots,\ol u_{N-1})\in U\big\}.$$
These are the set of choices of $u_{k}$ that are consistent with the
preceding choices $u_0,\ldots,u_{k-1}$, and are also consistent with feasibility. The last stage
corresponds to the complete solutions $u=(u_0,\ldots,u_{N-1})$, with cost 
$G(u)$; see Fig.\ 3.1.  All other transitions in this DP problem formulation have cost 0.
Of course here the number of states typically grows exponentially with $N$, but we can still apply the constrained rollout algorithm to the preceding DP formulation, using a suitable base heuristic, which will be applied only $N$ times. As an example, in Section 5 the multidimensional assignment problem will be transformed to the format of problem \gendiscrete, prior to the application of constrained rollout.  

The following  example describes a treatment by rollout of a classical 0-1 integer programming problem.

\xdef\examplefacility{\exampl}\examplnum\show{myexample}

\figrnum\show{myfigure}

\beginexample{\examplefacility\ (Facility Location)}We are given a  candidate set of $N$ locations, and we want to place in some of these locations a
``facility" that will serve the needs of $M$ ``clients." Each client $i=1,\ldots,M$ has a demand $d_i$ for services that may be satisfied at a location $k=0,\ldots,N-1$ at a cost $a_{ik}$ per unit. If a facility is placed at location $k$, it has capacity to serve demand up to a known level $c_k$.

\topinsert
\centerline{\hskip0pc\includegraphics[width=3.5in]{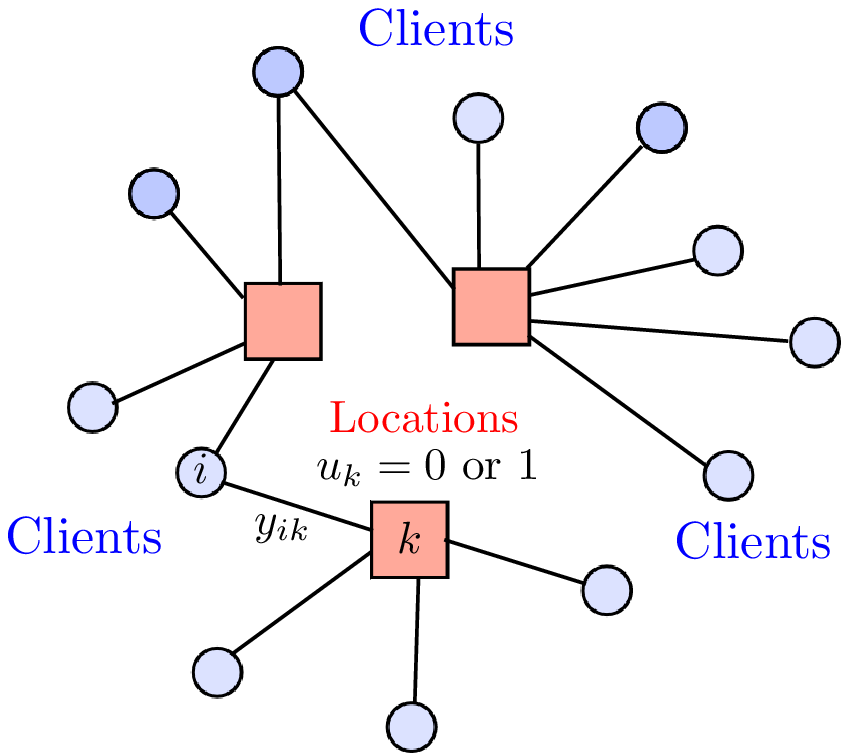}}
\vskip-1pc
\fig{0pc}{3.2.} {Schematic illustration of the facility location problem; cf.\ Example \examplefacility. Clients are matched to facilities, and the location of the facilities is subject to optimization.}\endinsert

We introduce a 0-1 integer
variable  
$u_k$ to indicate with $u_k=1$ that a facility is placed at location $k$ at a cost $b_k$ and with  $u_k=0$  that a facility is not placed at location $k$. Thus if $y_{ik}$ denotes the amount of demand of client $i$ to be served at facility $k$, the constraints are
$$\sum_{k=0}^{N-1} y_{ik}=d_i,\qquad i=1,\ldots,M,\xdef\constrone{\lab}\eqnum\show{oneo}$$
$$\sum_{i=1}^M y_{ik}\le c_ju_k,\qquad k=0,\ldots,N-1,\eqnum\show{oneo}$$
together with 
$$y_{ik}\ge0,\quad u_k\in\{0,1\},\qquad  i=1,\ldots,M,\ k=0,\ldots,N-1.\xdef\constrthree{\lab}\eqnum\show{oneo}$$
 We wish to minimize the cost
$$\sum_{i=1}^M\sum_{k=0}^{N-1} a_{ik}y_{ik}+\sum_{k=0}^{N-1} b_k u_k\xdef\locationcost{\lab}\eqnum\show{oneo}$$
subject to the preceding constraints. The essence of the problem is to place enough facilities at favorable locations to satisfy the  clients' demand at minimum cost. This can be a very difficult mixed integer programming problem.

On the other hand, when all the variables $u_k$ are fixed at some 0 or 1 values, the problem belongs to the class of {\it linear transportation} problems (see e.g., [Ber98]), and can be solved by fast polynomial algorithms. Thus the essential difficulty of the problem is how to select the sequence of variables $u_k$, $j=0,\ldots,N-1$. This can be viewed as a discrete optimization problem of the type shown in Fig.\ 3.2. In terms of the notation of this figure, the control components are $u_0,\ldots,u_{N-1}$, where $u_k$ can take the two values 0 or 1. 

To address the problem  by rollout, we must define a base heuristic at a ``state"
 $(u_0,\ldots,u_{k-1})$, where $u_{j}=1$ or $u_{j}=0$ specifies that a facility is or is not placed at location $j$, respectively. A suitable base heuristic at that state  is to place a facility at all of the remaining locations (i.e., $u_j=1$ for $j=k+1,\ldots,N-1$), and its cost is obtained by solving the corresponding linear transportation problem of minimizing the cost \locationcost\ subject to the constraints \constrone-\constrthree, with the variables $u_j$, $j=0,\ldots,k-1$, fixed at the previously chosen values, and the variables $u_j$, $j=k,\ldots,N$, fixed at 1.
 
To illustrate, at the initial state where no placement decision has been made, we set $u_0=1$ (a facility is placed at location 0) or $u_0=0$ (a facility is not placed at location 0), we solve the  two corresponding transportation problems, and we fix $u_0$, depending on which of the two resulting costs is smallest. Having fixed the status of location 0, we repeat with location 1, set the variable $u_1$ to 1 and to 0, solve the corresponding two transportation problems, and fix $u_1$, depending on which of the two resulting costs is smallest, etc.

It is easily seen that if the initial base heuristic choice (placing a facility at every candidate location) yields a feasible solution, i.e., 
$$\sum_{i=1}^Md_i\le \sum_{k=0}^{N-1} c_k,$$
the rollout algorithm will yield a feasible solution with cost that is no larger than the cost corresponding to the initial application of the base heuristic. In fact it can be verified that the base heuristic here is sequentially improving, so it is not necessary to use the fortified version of the algorithm. Regarding computational costs, the number of transportation problems to be solved is at first count $2N$, but it can be reduced to $N+1$ by exploiting the fact that one of the two transportation problems at each stage after the first has been solved at an earlier stage. It is finally worth noting, for readers that are familiar with the integer programming method of branch-and-bound, that the graph of  Fig.\ 3.1 corresponds to the branch-and-bound tree for the problem, so the rollout algorithm amounts to a quick (and imperfect) method to traverse the branch-and-bound tree. This observation may be useful if we wish to use integer programming techniques to add improvements to the rollout algorithm.

We finally note that the rollout algorithm requires the solution of many linear transportation problems, with fairly similar data. It is thus important to use an algorithm that is capable of using effectively the final solution of one transportation problem as a starting point for the solution of the next. The auction algorithm for transportation problems (Bertsekas and Casta\~non [Ber89]) is particularly well-suited for this purpose.
\endexample

\vskip-1.5pc

\section{Constrained Multiagent Rollout}
\mark{Constrained Multiagent Rollout}
\vskip-0.5pc

\pn Let us  assume a special structure of the control space, where the control $u_k$ consists of $m$ components, $u_k=(u_k^1,\ldots,u_k^m)$, each belonging to a corresponding set $U_k^\ell$, $\ell=1,\ldots,m$. 
Thus  the control space at stage $k$ is the Cartesian product
$$U_k=U_k^1\times \cdots \times U_k^m.\xdef\eqthree{\lab}\eqnum\show{oneo}$$
We refer to this as the {\it multiagent case\/}, motivated by the special case where each component $u_k^\ell$, $\ell=1,\ldots,m$, is chosen by a separate agent $\ell$ at stage $k$. Then the rollout minimization \mincontrt\ involves the computation and comparison of as many as $n^m$ terms $G\big(T_k(\tl y_k,u_k)\big)$, where $n$ is the maximum number of elements of the sets $U_k^\ell$ [so that $n^m$ is an upper bound to the number of controls in the control space $U_k$, in view of its Cartesian product structure \eqthree]. Thus the  rollout algorithm requires order $O(n^m)$ applications of the base heuristic per stage.

In this section we construct an alternative rollout algorithm that  achieves the cost improvement property \costimprove\ with much smaller computational cost, namely order $O(nm)$ applications of the base heuristic per stage. A key idea here is that the computational requirements of the rollout one-step minimization \mincontrt\ are proportional to the number of controls and are independent of the size of the state space. This motivates a reformulation of the problem, first suggested in the neuro-dynamic programming book [BeT96], Section 6.1.4, whereby control space complexity is traded off with state space complexity, by ``unfolding" the control $u_k$ into its $m$ components, which are applied one {\it agent-at-a-time} rather than {\it all-agents-at-once\/}. We describe this idea next within the context of this paper; see [Ber19b] for a related discussion, which also considers stochastic multiagent problems and associated rollout algorithms with a cost improvement guarantee over the base heuristic.

\subsection{Trading off Control Space Complexity with State Space 
Complexity}
\pn We noted that a major issue in rollout is the  minimization over $u_k$
$$\tl u_k\in\arg\min_{u_k\in U_k(\tl y_k)}G\big(T_k(\tl y_k,u_k)\big),$$ 
[cf.\ Eq.\ \mincontrt], which can be very time-consuming when the size of the control space is large. In particular, in the multiagent case when $u_k=(u_k^1,\ldots,u_k^m),$ the time to perform this minimization is typically exponential in $m$. To deal with this, we reformulate the problem by breaking down the collective decision $u_k$ into $m$ individual component decisions, thereby reducing the complexity of the control space while
increasing the complexity of
the state space. The potential advantage is that the extra state space 
complexity does not affect the computational requirements of rollout.

\xdef\figunfolded{\figr}\figrnum\show{myfigure}

To this end, we introduce a modified but equivalent problem, involving one-at-a-time agent control selection. In particular, at the generic state $x_k$, we break down the control $u_k$ into the 
sequence
of the $m$ controls $u_k^1,u_k^2,\ldots,u_k^m$, and between $x_k$ and the next state $x_{k+1}=f_k(x_k,u_k)$, we introduce artificial 
intermediate ``states" $(x_k,u_k^1),(x_k,u_k^1,u_k^2),\ldots,(x_k,u_k^1,\ldots,u_k^{m-1})$, and corresponding 
transitions. The choice of the last control component $u_k^m$ 
at ``state" $(x_k,u_k^1,\ldots,u_k^{m-1})$ marks the transition at cost $g_k(x_k,u_k)$ to the next state $x_{k+1}=f_k(x_k,u_k)$ according to 
the system equation; see Fig.\ \figunfolded.

\topinsert
\centerline{\hskip0pc\includegraphics[width=5in]{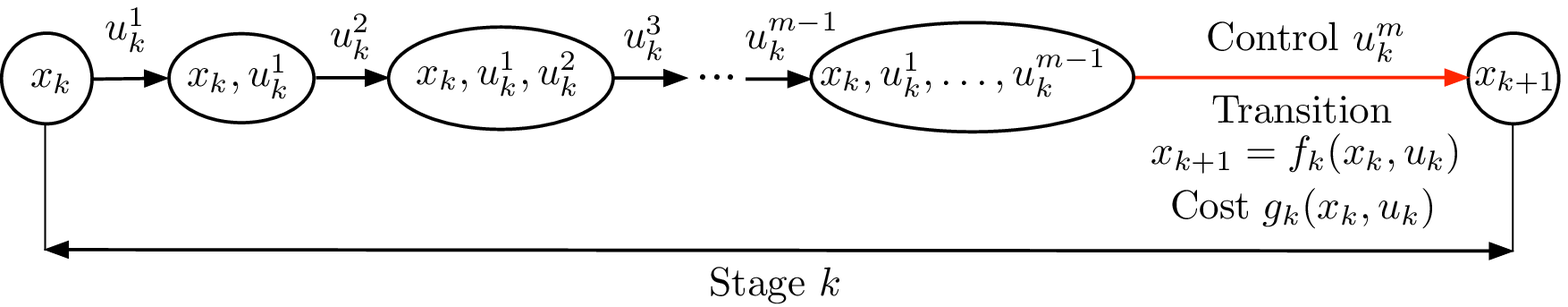}}
\vskip-1pc
\hskip-2pc\fig{0pc}{4.1.} {Equivalent formulation of the $N$-stage optimal control problem for the case where the control $u_k$ consists of $m$ components $u_k^1,u_k^2,\ldots,u_k^m$:
$$u_k=(u_k^1,\ldots,u_k^m)\in U_k^1(x_k)\times \cdots \times U_k^m(x_k).$$
The figure depicts the $k$th stage transitions. Starting from state $x_k$, we generate the intermediate states $(x_k,u_k^1),(x_k,u_k^1,u_k^2),\ldots,(x_k,u_k^1,\ldots,u_k^{m-1})$, using the respective controls $u_k^1,\ldots, u_k^{m-1}$. The final control $u^m$ leads from $(x_k,u_k^1,\ldots,u_k^{m-1})$ to $x_{k+1}=f_k(x_k,u_k)$, and a cost $g_k(x_k,u_k)$ is incurred.}
\endinsert

It is evident that this reformulated problem is equivalent to the original, since any control choice that is possible in one problem is also possible in the other problem, while the cost structures of the two problems are essentially the same. The motivation for the reformulated problem is that the control space is simplified 
at the expense of
introducing $m-1$ additional layers of states.  However, the key point is that {\it the increase in size of the state 
space does not adversely affect the operation of rollout\/}. 

\subsection{Multiagent Rollout Algorithm}

\pn Consider now the constrained rollout algorithm of Section 2 applied to the reformulated problem shown in Fig.\ \figunfolded, with a base heuristic suitably modified so that it completes a partial trajectory of the form
$$\big(y_k, (x_k,u_k^1),(x_k,u_k^1,u_k^2),\ldots,(x_k,u_k^1,\ldots,u_k^{\ell})\big),\qquad \ell=1,\ldots,m.$$
The algorithm involves a minimization over only one control component  at the state $x_k$ and at each of the intermediate states 
$$(x_k,u_k^1),(x_k,u_k^1,u_k^2),\ldots,(x_k,u_k^1,\ldots,u_k^{m-1}).$$
 In particular, {\it for each stage $k$, the algorithm requires a sequence of $m$ minimizations, one over each of the  control components $u_k^1,\ldots,u_k^m$, with the past controls already determined by the rollout algorithm, and the future controls determined by running the base heuristic.\/}  Assuming a maximum of $n$ elements in the control component spaces $U_k^\ell$, $\ell=1,\ldots,m$, the computation required at each stage $k$ is of order $O(n)$ for each of the ``states" 
$$x_k,(x_k,u_k^1),\ldots,(x_k,u_k^1,\ldots,u_k^{m-1}),$$
 for a total of order $O(nm)$ computation. 

To elaborate, for all $k$ and $\ell\le m$ at the current partial trajectory
$$(\tl x_0,\tl u_0,\ldots,\tl x_k,\tl u_k^1,\ldots,\tl u_k^{\ell-1}),$$ 
and for each of the controls $u_k^\ell$, we use the base heuristic to generate a complementary partial trajectory
$$(u_k^{\ell+1},\ldots,u_k^m,x^{k+1},u^{k+1},\ldots,x^{N-1},u^{N-1},x_N),\xdef\completraj{\lab}\eqnum\show{oneo}$$
 up to stage $N$. We then select the control $\tl u_k^\ell$ for which the resulting complete trajectory is feasible and has minimum cost. There is also a fortified version of this algorithm, which is similar to the one described in Section 3.

Note that the base heuristic used in the reformulated problem must be capable of generating a complementary partial trajectory of the form \completraj, starting from any partial trajectory of states and controls. 
Note also that instead of selecting the agent controls in a fixed order, it is
possible to change the order at each stage $k$. In fact it is possible to optimize over multiple orders at the same stage. 

\vskip-1.pc

\section{Application to Multidimensional Assignment}
\mark{Application to Multidimensional Assignment}

\vskip-1pc

\pn In this section we demonstrate the application of constrained multiagent rollout within the context of the classical multidimensional assignment problem. This problem is representative of a broad class of layered graph problems, which involve both a temporal and a spacial allocation structure, so that both the dynamic system character and the multiagent character of our algorithms come into play.  

Multidimensional assignment problems involve graphs consisting of $N+1$ subsets of nodes ($N\ge 2$), denoted ${\cal N}_0,{\cal N}_1,\ldots,{\cal N}_N$, and referred to as {\it layers\/}. The arcs of the graphs are directed and are of the form $(i,j)$, where $i$ is a node in a layer ${\cal N}_k$, $k=0,1,\ldots,N-1$, and $j$ is a node in the corresponding next layer ${\cal N}_{k+1}$. Thus we have a directed graph whose nodes are arranged in $N+1$ layers and the arcs connect the nodes of each layer to the nodes in their adjacent layers; see Fig.\ 5.1. Here for simplicity, we assume that each of the layers ${\cal N}_k$ contains $m$ nodes, and that there is a unique arc connecting each node in a given layer with each of the nodes of the adjacent layers.

\topinsert
\centerline{\hskip0pc\includegraphics[width=4in]{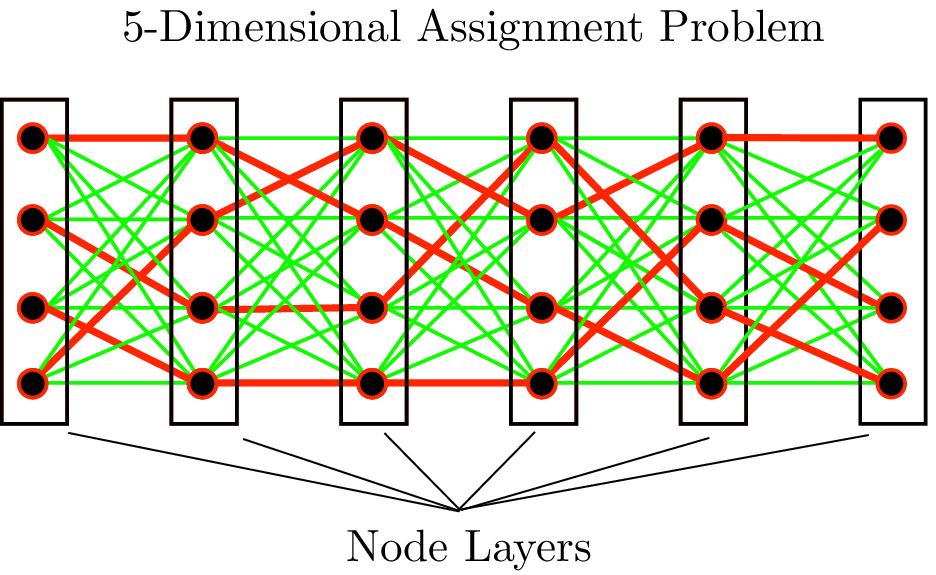}}
\vskip-1pc
\hskip-2pc\fig{0pc}{5.1.} {Illustration of the graph of an $(N+1)$-dimensional assignment problem (here $N=5$). There are $N+1$ node layers each consisting of $m$ nodes (here $m=4$). Each grouping consists of $N+1$ nodes, one from each layer, and $N$ corresponding arcs. An $(N+1)$-dimensional assignment consists of $m$ node-disjoint groupings, where  each node belongs to one and only one grouping (illustrated in the figure with thick red lines). For each grouping, there is an  associated cost, which depends on the $(N+1)$-tuple of arcs that comprise the grouping. The cost of an $(N+1)$-dimensional assignment is the sum of the costs of its $m$ groupings. In the separable case, the cost of a grouping separates into the sum of its $N$ arc costs, and the problem can be solved by solving $N$ decoupled 2-dimensional assignment problems.}
\endinsert

We consider subsets of nodes, referred to as {\it groupings\/}, which contain a single node from every layer, and we assume that every grouping  is associated with a given cost. A partition of the set of nodes into $m$ disjoint groupings, so that each node belongs to one and only one grouping is called an {\it $(N+1)$-dimensional assignment\/}. For each grouping, there is an  associated cost, which depends on the $N$-tuple of arcs that comprise the grouping. The cost of an $(N+1)$-dimensional assignment is the sum of the costs of its $m$ groupings. The problem is to find an $(N+1)$-dimensional assignment of minimum cost.

This is a difficult combinatorial problem with many applications.  
An important special case arises in
the context of multi-target tracking and data association; see Blackman [Bla86],
Bar-Shalom and Fortman [BaF88], Bar-Shalom [Bar90], Pattipati, Deb, Bar-Shalom, and Washburn [PDB92],
Poore [Poo94], Poore and Robertson [PoR97], Popp, Pattipati, and Bar-Shalom [PPB01], and Choi, Brunet, and How [CBH09]. Other challenging combinatorial problems, such as multi-vehicle routing and multi-machine scheduling problems, also share the spatio-temporal type of structure, and are thus well suited for the application of our constrained multiagent rollout approach. 
For more discussion of related combinatorial applications, we refer to Chapter 10 of the author's network optimization book [Ber98]. Generally, the fine details of such a problem will determine the choice of a suitable base heuristic. 

We note that there are several variants of the multidimensional assignment problem illustrated in Fig.\ 5.1, which are well-suited for the application of constrained rollout. For example, these variants may involve unequal numbers of nodes in each layer, or a sparse structure where some of the possible arcs connecting nodes of adjacent layers are missing. Moreover, there may be cost coupling between collections of groupings that depends on the groupings' compositions. In this paper we will focus on the case where the layers have equal numbers of nodes and where the cost of each grouping depends exclusively on the $N+1$  nodes that comprise the grouping. This structure favors the use of base heuristics that rely  on solution of 2-dimensional assignment problems.

\subsection{Three-Dimensional Assignment}

\pn To simplify the presentation, we will first focus on the 3-dimensional assignment special case ($N=2$), and for descriptive purposes, we will associate the nodes of the three layers with ``jobs," ``machines," and ``workers," respectively. Suppose that the performance of a job $j$ requires a single machine $\ell$ and a single worker $w$ (which cannot be shared by any other job),
and that there is a given cost
$a_{j\ell w}$ corresponding to the triplet
$(j,\ell,w)$. Given a set of $m$ jobs, a set of $m$ machines, and a set of $m$ workers, we want to find a collection of $m$ job-machine-worker triplets that has
minimum total cost. This problem is quite challenging, and is well-suited for demonstration of our constrained rollout approach as it has both a temporal and a spacial character.

To transcribe the problem to the optimal control format of this paper, we use the formulation of Fig.\ 3.1, and assume that there is only one state at each of three stages [the respective collections of jobs (for the first stage), machines (for the second stage), and workers (for the third stage)], and two decisions to make (the assignment of jobs to machines and the assignment of machines to workers); see Fig.\ 5.2. Each of the decisions consists of $m$ components, the $m$ outgoing arcs  from the $m$ nodes corresponding to the stage. Thus, the application of the multiagent rollout algorithm of Section 4 will involve two stages, a state space consisting of a single element for each state, and a control at each stage that consists of $m$ components (the choice of machine to assign to each job in the first stage, and the choice of worker to assign to each machine in the second stage). These components are computed in sequence according to some predetermined order, which without loss of generality we will assume to be the natural order $1,2,\ldots,m$.

\topinsert
\centerline{\hskip0pc\includegraphics[width=3.5in]{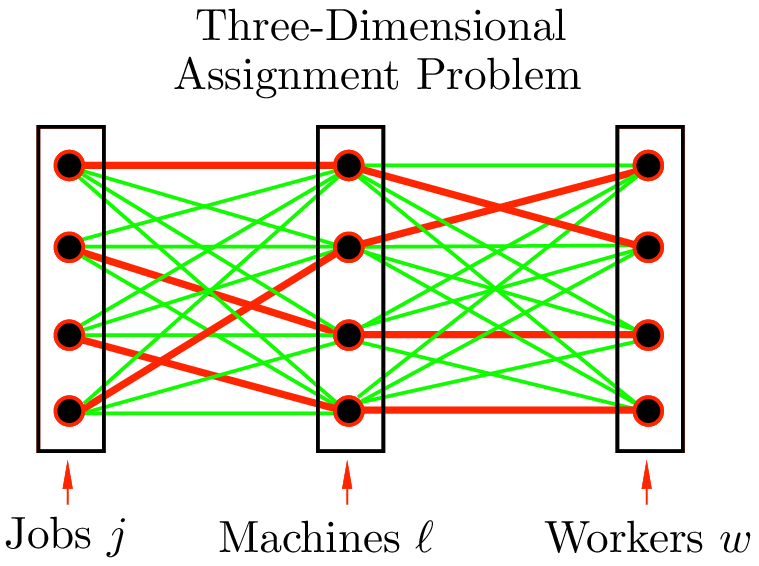}}
\vskip-1pc
\hskip-2pc\fig{0pc}{5.2.} {A 3-dimensional assignment problem consisting of assigning each job $j$ to a machine $\ell$ and to a worker $w$ at cost $a_{j\ell w}$. Each machine assigned to exactly one job and exactly one worker. Here $N=2$ and $m=4$.}
\endinsert

An important and particularly favorable special case of the problem arises when the costs $a_{j\ell w}$
have the {\it separable} form
$$ a_{j\ell w}=\b_{j\ell}+\g_{\ell w},$$
where $\b_{j\ell}$ and $\g_{\ell w}$ are given scalars. In this case, there is no coupling between jobs and
workers, and the problem can be efficiently (polynomially) solved by solving two decoupled (2-dimensional) assignment 
problems: one involving the pairing of jobs and machines, with the $\b_{j\ell}$ as costs, and the
other involving the pairing of machines and workers, with the $\g_{\ell w}$ as costs. In general,
however, the 3-dimensional assignment problem is a difficult integer programming problem, 
for which there is no known polynomial algorithm.

\subsection{The Enforced Separation Heuristic and Three-Dimensional Assignment}

\pn The  separable case motivates a simple heuristic approach for the nonseparable 3-dimensional assignment problem, which is well-suited for the use of constrained rollout. The base heuristic is a two-stage procedure that is based on relaxing each the grouping constraints, by first focusing on assigning machines to workers, and then focusing on assigning jobs to machines, after suitably modifying the costs $a_{j\ell w}$ to make them separable. In
particular, we first relax the constraints on the jobs by assuming that each machine-worker pair can be assigned to the most favorable job. Mathematically, this relaxed
problem takes the polynomially solvable 2-dimensional assignment form
$$\eqalignno{
\hbox{\rm minimize}\quad &\sum_{w=1}^m\sum_{\ell=1}^mc_{\ell w}v_{\ell w}&\cr
\hbox{\rm subject to\ \ }
&\sum_{\ell=1}^m v_{\ell w}=1,\qquad \forall\ w=1,\ldots,m,&\cr
&\sum_{w=1}^mv_{\ell w}=1,\qquad \forall\ \ell=1,\ldots,m,&\cr
&v_{\ell w}=0\hbox{ or }1,\qquad \forall\ \ell,\ w=1,\ldots,m,&\cr
}$$
where
$$c_{\ell w}=\min_{j=1,\ldots,m} a_{j\ell w},\xdef\mincosts{\lab}\eqnum\show{spconst}$$
 $v_{\ell w}$ are the variables of the problem, and $v_{\ell w}=1$ indicates that machine $\ell$ is assigned to worker $w$, so the constraints of the above problem enforce the condition that each machine is assigned to one and only one worker.\footnote{\dag}{An alternative is to define each cost $c_{\ell w}$ as a ``representative" cost $a_{j_{\ell w}\ell w}$ (for some specially selected job $j_{\ell w}$). Such an alternative may become attractive when extensions of enforced separation are considered for $(N+1)$-dimensional problems with large $N$. Then the analog of Eq.\ \mincosts\ will involve minimization over the  (exponential in $N$) number of all $(N-2)$-tuples of graph arcs that can form a grouping with $(\ell,w)$, and may become very costly.}

For each $w$, let $w_\ell$ be
the worker assigned to machine $\ell$, according to the solution of this problem. We can now optimally
assign jobs $j$ to machine-worker pairs $(\ell,w_\ell)$ by using as assignment costs 
$$b_{j\ell}=a_{j\ell w_\ell},$$
and obtain a (suboptimal) 3-dimensional assignment $\{(j_\ell,\ell,w_\ell)\mid \ell=1,\ldots,m\}$. It can be seen that this
approach amounts to {\it enforced separation\/}, whereby we replace the costs $a_{j\ell w}$ with
the separable approximations
$b_{j\ell}+c_{\ell w}$. In fact, it can be shown that if the problem is {\it $\e$-separable\/}, 
in the sense
that for some (possibly unknown) $\ol \b_{j\ell}$ and $\ol \g_{\ell w}$, and some $\e\ge0$, we have
$$|\ol \b_{j\ell}+\ol \g_{\ell w}-a_{j\ell w}|\le\e,\qquad\forall\ j,\,\ell,\,w=1,\ldots,m,$$
then the assignment $\{(j_\ell,\ell,w_\ell)\mid \ell=1,\ldots,m\}$ obtained using the preceding enforced
separation approach achieves the optimal cost of the problem within $4m\e$ (see the author's network optimization book [Ber98], Exercise 10.31).

The enforced separation approach is simple and can be generalized to problems with more than two stages as we will discuss shortly. 
Moreover, enforced separation heuristics also apply to several variants of the
multidimensional assignment problem.   For example, we may have transportation-type constraints, where
multiple jobs can be performed on the same machine, and/or multiple
machines can be operated by a single worker. In this case, our preceding discussion of the
enforced separation heuristic applies similarly, except that we need to solve 2-dimensional transportation problems rather than 2-dimensional assignment problems. 

\subsection{Using Enforced Separation as a Base Heuristic}

\pn We will now describe the use of enforced separation as a base heuristic in the context of constrained rollout. The 3-dimensional assignment problem is posed as an optimal control problem involving $m+1$ sequential choices: the machines assigned to the jobs are first selected one-by-one in some fixed order ($m$ sequential choices), and then the workers assigned to the machines are selected simultaneously. To connect with our earlier optimal control formulation, trajectories here consist of an artificial initial state, the $m$ successive choices of job-machine pairs (these correspond to the controls $u_0^1,\ldots,u_0^m$), and then the $m$-tuple of machine-worker pairs (these comprise the control $u_1$).
For each of the first $m$ choices a job is selected and the machine to be assigned to this job is fixed by the rollout algorithm, through the use of the base heuristic of enforced separation. At the last stage the $m$ machines are assigned simultaneously to workers using a 2-dimensional assignment algorithm.

To illustrate the fortified rollout algorithm, at the artificial initial state where no job-machine or machine-worker pair has been fixed, the enforced separation heuristic as described above is used to generate a (suboptimal) initial 3-dimensional assignment, which serves as the initial tentative best trajectory, and has  cost denoted by $\hat S$. 

In the first $m$ rollout stages we select  in sequence each job $j=1,\ldots,m$, and we select a machine $\ell$ to assign to it, by using the enforced separation heuristic. The first rollout stage is as follows:

\nitem{Stage 1.1:} We take job $j=1$, and fix its assignment to machine 1. We then apply the enforced separation heuristic by solving two 2-dimensional assignment problems. The first of these involves the assignment of machines to workers using as costs 
$$c_{\ell w}=\cases{a_{11w}&if $\ell=1$,\cr
\min_{j=2,\ldots,m} a_{j\ell w}&if $\ell\ne1$,\cr}$$
[cf.\ Eq.\ \mincosts]. We thus obtain an assignment of machines to workers of the form $(\ell,w_\ell)$, $\ell=1,\ldots,m$. Having fixed the workers to be assigned to machines, we solve the 2-dimensional assignment problem of assigning the jobs $2,\ldots,m$ to the machines $2,\ldots,m$, where the costs are 
$$b_{j\ell}=a_{j\ell w_\ell},\qquad j,\ell=2,\ldots,m.$$
By ``joining" the solutions of the two 2-dimensional problems just described, we obtain a 3-dimensional assignment that consists of $m$ job-machine-worker groupings whose cost we call $S_1$.

\nitem{Stage 1.2:} We take job $j=1$, and fix its assignment to machine 2. We then apply the enforced separation heuristic by solving two 2-dimensional assignment problems. The first involves the assignment of machines to workers using as costs 
$$c_{\ell w}=\cases{a_{12w}&if $\ell=2$,\cr
\min_{j=2,\ldots,m} a_{j\ell w}&if $\ell\ne 2$,\cr}$$
[cf.\ Eq.\ \mincosts]. We obtain an assignment of machines to workers of the form $(\ell,w_\ell)$, $\ell=1,\ldots,w$. We then solve the 2-dimensional assignment problem of assigning the jobs $2,3,\ldots,m$ to the machines $1,3,\ldots,m$, where the costs are 
$$b_{j\ell}=a_{j\ell w_\ell},\qquad j=2,3,\ldots,m,\ \ell=1,3,\ldots,m.$$
By ``joining" the solutions of the two 2-dimensional problems, we obtain a 3-dimensional assignment, whose cost we call $S_2$.

\nitem{Stage 1.t:} For $t=3,\ldots,m$, we continue the process described above, where we fix the assignment of job 1 to machine $t$. We then apply the enforced separation heuristic by solving two 2-dimensional assignment problems, similar to the ones above: first assigning machines to workers using costs
$$c_{\ell w}=\cases{a_{1tw}&if $\ell=t$,\cr
\min_{j=2,\ldots,m} a_{j\ell w}&if $\ell\ne t$,\cr}$$
and obtaining an assignment of machines to workers of the form $(\ell,w_\ell)$, $\ell=1,\ldots,w$. We then solve the 2-dimensional assignment problem of assigning the jobs $2,\ldots,m$ to the machines $1,\ldots,t-1,t+1,\ldots,m$, where the costs are 
$$b_{j\ell}=a_{j\ell w_\ell},\qquad j=2,3,\ldots,m,\ \ell=1,\ldots,t-1,t+1,\ldots,m.$$
By ``joining" the solutions of the two 2-dimensional problems, we obtain a 3-dimensional assignment, whose cost we call $S_t$.
\smskip

\nitem{Stage 1:} We now have $m$ 3-dimensional assignments, with corresponding costs $S_1,\ldots,S_m$, where job 1 is fixed to the machines $1,\ldots,m$, respectively. We choose the machine $\tl \ell$ for which $S_{\tl\ell}$ is minimized over $\ell=1,\ldots$, and we permanently assign job 1 to machine $\tl \ell$ if $S_{\tl\ell}$ is less or equal to $\hat S$, the cost of the enforced separation heuristic applied to the artificial initial state (the current  tentative best trajectory), and we also adopt the corresponding 3-dimensional assignment as the new tentative best trajectory. Otherwise, we set the assignment of job 1 to a machine according to the current  tentative best trajectory, which we leave unchanged.

The preceding procedure, the first step of constrained rollout, required the solution of $2m$ 2-dimensional assignment problems, and yielded a permanent assignment of job 1 to a machine. The procedure is then repeated for job 2, taking into account that the assignment of job 1 to a machine has been fixed. This requires similarly the solution of $2(m-1)$ 2-dimensional assignment problems, and yields a permanent assignment of job 2 to a machine, and an update of the current tentative trajectory. Repeating the procedure with jobs $3,\ldots,m$ in sequence, we obtain a permanent   assignment of all the jobs to machines, and the corresponding 3-dimensional assignment (which has minimum cost over all the 3-dimensional assignments generated, in view of the use of fortified rollout). The total number of 2-dimensional assignment problems thus solved is
$$2m+2(m-1)+2(m-2)+\cdots+2=m^2.$$
Finally, given the permanent assignment of all the jobs to machines, say $(j_\ell,\ell)$, $\ell=1,\ldots,m$, we obtain the (permanent) assignment of workers to job-machine pairs, by using as costs the scalars
$$c_{\ell w}=a_{j_\ell \ell w},\qquad \ell,w=1,\ldots,m.$$
At this point we have obtained by rollout the final suboptimal 3-dimensional assignment, which by construction has the cost improvement property: it has no larger cost than the one obtained by the enforced separation base heuristic starting from the artificial initial condition.

Thus the total number of 2-dimensional assignment problems to be solved by the rollout algorithm is $m^2+1$. Each of these problems can be solved very fast using any one of a number of methods. However, because these problems and their solutions are similar, it is important to use a method that can exploit this similarity. A particularly favorable method in this regard is the author's auction algorithm [Ber79] (see the book [Ber98] for a detailed development; we provide an introductory review of auction algorithms in the appendix to this paper).   The auction algorithm uses a price variable for each node, such as a worker or a machine, and then adjusts the prices through an auction-like process to achieve a form of economic equilibrium. One can then use the final prices obtained for one 2-dimensional assignment problem as an efficient starting point for the solution of a related 2-dimensional assignment problem.  For this reason, the auction algorithm and its variations have been widely adopted for use in solving 2-dimensional assignment problems in the context of multidimensional assignment algorithms used for multitarget tracking applications, among others (see [PDB92], [Poo94], [PoR97], [PPB01], [CBH09]). 

\subsection{Enforced Separation and Constrained Rollout for Multidimensional Assignment}

\pn Let us now consider briefly the extension of the constrained rollout algorithm just described to 
$(N+1)$-dimensional assignment problems with $N>2$. Here we will use an extension of the 3-dimensional enforced separation heuristic. We start again from the last stage, solve the last 2-dimensional assignment problem of the last stage by modifying the arc costs according to the analog of the minimization formula 
$$c_{\ell w}=\min_{j=1,\ldots,m} a_{j\ell w},$$
[cf.\ Eq.\ \mincosts]. The only difference is that instead of minimizing $a_{j\ell w}$ over jobs $j$ as above, we minimize over all $(N-1)$-tuples of nodes of the groupings whose final two nodes are $(\ell,w)$. Once the assignments of the last stage are fixed, a similar procedure can be used to fix the  assignments of the next-to-last stage, and so on. The total number of 2-dimensional assignment problems to be solved is $(m+1)(N-2)$. Thus the base heuristic's computation time is polynomial in both $m$ and $N$.

The enforced separation heuristic just described for the artificial initial condition, can be used in suitably modified form for rollout, with the assignments of some stages fixed permanently by rollout, and additional assignments fixed one by one, by applying the enforced separation heuristic, with suitably modified costs that take into account the already fixed assignments.

\vskip-1pc

\section{Concluding Remarks}
\mark{Concluding Remarks}

\vskip-1pc

\pn We have proposed a rollout algorithm for constrained deterministic DP, which is well suited for the suboptimal solution of challenging discrete optimization problems. Under suitable assumptions, we have shown a cost improvement property: the rollout algorithm produces a feasible solution, whose cost is no worse than the cost of the solution produced by the base heuristic. 

We have also proposed an efficient variant of the algorithm for multiagent problems, where the control at each stage consists of multiple components. We have shown that the cost improvement property is preserved, and we have applied the algorithm to layered graph problems that involve both a spatial and a temporal structure. In particular, we have focused on multidimensional assignment, using the auction algorithm for 2-dimensional assignment as a base heuristic.

\vskip-1pc

\def\ref{\vskip1pt\pn}

\section{References}
\mark{References}
\vskip-1pc

\ref[AnH14] Antunes, D., and Heemels, W.P.M.H., 2014.\ ``Rollout Event-Triggered Control: Beyond Periodic Control Performance," IEEE Transactions on Automatic Control, Vol.\ 59, pp.\ 3296-3311.

\ref[BBG13] Bertazzi, L., Bosco, A., Guerriero, F., and Lagana, D., 2013.\ ``A Stochastic Inventory Routing Problem with Stock-Out," Transportation Research, Part C, Vol.\ 27, pp.\ 89-107.

\ref [BCE95] Bertsekas, D.\ P., Casta\~ non, D.\ A., Eckstein, J., and Zenios, S.,
1995.\ ``Parallel Computing in Network Optimization," Handbooks in OR and
MS,  Ball, M.\ O., Magnanti, T.\ L., Monma, C.\ L., and Nemhauser, G.\ L.\
(eds.),  Vol.\ 7, North-Holland, Amsterdam, pp.\ 331-399.

\ref [BTW97] Bertsekas, D.\ P., Tsitsiklis, J.\ N., and Wu, C., 1997.\ ``Rollout Algorithms for
Combinatorial Optimization,'' Heuristics, Vol.\ 3, pp.\ 245-262.

\ref [BaF88] Bar-Shalom, Y., and Fortman, T.\ E., 1988.\ Tracking and Data
Association, Academic Press, N.\ Y.

\ref [BaF90] Bar-Shalom, Y., 1990.\ Multitarget-Multisensor Tracking: Advanced Applications, Artech House, Norwood, MA.

\ref [BeC89] Bertsekas, D.\ P., and Casta\~ non, D.\ A., 1989.\ ``The Auction
Algorithm for Transportation Problems,'' Annals of Operations Research,
Vol.\ 20, pp.\ 67-96.

\ref [BeC91] Bertsekas, D.\ P., and Casta\~ non, 1991.\ ``Parallel
Synchronous and Asynchronous Implementations of the Auction Algorithm,'' 
Parallel Computing, Vol.\ 17, pp.\ 707-732.

\ref[BeC99] Bertsekas, D.\ P., and  Casta\~ non, D.\ A., 1999.\ ``Rollout Algorithms for
Stochastic Scheduling Problems," Heuristics, Vol.\ 5, pp.\ 89-108. 

\ref[BeC08] Besse, C., and Chaib-draa, B., 2008.\ ``Parallel Rollout for Online Solution of DEC-POMDPs," Proc.\ of 21st International FLAIRS Conference, pp. 619-624.

\ref[BeG97] Beraldi, P., and Guerriero, F., 1997.\ ``A Parallel Asynchronous Implementation of the
Epsilon-Relaxation Method for the Linear Minimum Cost Flow Problem,"
Parallel Computing, Vol.\ 23, pp.\ 1021-1044.

\ref[BeL14] Beyme, S., and Leung, C., 2014.\ ``Rollout Algorithm for Target Search in a Wireless Sensor Network," 80th Vehicular Technology Conference (VTC2014), IEEE, pp.\ 1-5.

\ref[BeS18] Bertazzi, L., and Secomandi, N., 2018.\ ``Faster Rollout Search for the Vehicle Routing Problem with Stochastic Demands and Restocking," European J.\ of Operational Research, Vol.\ 270, pp.\  487-497.

\ref [BeT89]
 Bertsekas, D.\ P., and Tsitsiklis, J.\ N., 1989.\ Parallel and
Distributed Computation: Numerical Methods, Prentice-Hall, Englewood Cliffs,
N.\ J.\ (republished in 1997 by Athena Scientific, Belmont, MA).

\ref [BeT96]  Bertsekas, D.\ P., and Tsitsiklis, J.\ N., 1996.\ Neuro-Dynamic Programming,
Athena Scientific, Belmont, MA.

\ref[Ber79] Bertsekas, D.\ P., 1979.\ ``A Distributed Algorithm for the Assignment Problem," Lab.\ for Information and Decision Systems Report, MIT, May 1979.

\ref [Ber86] Bertsekas, D.\ P., 1986.\ ``Distributed Relaxation Methods for
Linear Network Flow Problems,'' Proceedings of 25th IEEE Conference on Decision
and Control, Athens, Greece, pp.\ 2101-2106.

\ref [Ber88] Bertsekas, D.\ P., 1988.\ ``The Auction Algorithm: A Distributed
Relaxation Method for the Assignment Problem,'' Annals of Operations
Research, Vol.\ 14, pp.\ 105-123.

\ref [Ber91] Bertsekas, D.\ P., 1991.\ Linear Network Optimization: Algorithms
and Codes, MIT Press, Cambridge, MA.

\ref [Ber92] Bertsekas, D.\ P., 1992.\ ``Auction Algorithms for Network Flow
Problems: A Tutorial Introduction," Computational Optimization and
Applications, Vol.\ 1, pp.\ 7-66.

\ref [Ber93] Bertsekas, D.\ P., 1993.\ ``Mathematical Equivalence of the Auction Algorithm for
Assignment and the $\e$-Relaxation (Preflow-Push) Method for Min Cost Flow," in Large Scale
Optimization: State of the Art, Hager, W.\ W., Hearn, D.\ W., and Pardalos, P.\ M.\ (eds.),
Kluwer, Boston, pp.\ 27-46.

\ref[Ber97] Bertsekas, D.\ P., 1997.\ ``Differential Training of Rollout Policies,"
Proc. of the 35th Allerton Conference on Communication, Control, and Computing,
Allerton Park, Ill.

\ref [Ber98] Bertsekas, D.\ P., 1998.\ Network Optimization: Continuous and Discrete Models, Athena
Scientific, Belmont, MA (available in .pdf form through the author's web site).

\ref [Ber01] Bertsekas, D.\ P., 2001.\ ``Auction Algorithms," Encyclopedia of Optimization, Kluwer.

\ref[Ber05a] Bertsekas, D.\ P., 2005.\ ``Rollout Algorithms for Constrained Dynamic Programming," LIDS Report 2646, MIT.

\ref[Ber05b] Bertsekas, D.\ P., 2005.\ ``Dynamic Programming and Suboptimal Control: A Survey from ADP to MPC," European J.\ of Control, Vol.\ 11, pp.\ 310-334.

\ref [Ber17] Bertsekas, D.\ P., 2017.\ Dynamic Programming and Optimal Control, 4th Edition, Athena
Scientific, Belmont, MA.

\ref[Ber19a] Bertsekas, D.\ P., 2019. Reinforcement Learning and Optimal Control, Athena Scientific, Belmont, MA.

\ref[Ber19b] Bertsekas, D.\ P., 2019. ``Multiagent Rollout Algorithms," arXiv preprint arXiv:1910.00120.

\ref[Ber19c] Bertsekas, D.\ P., 2019. ``Robust Shortest Path Planning and Semicontractive Dynamic Programming," Naval Research Logistics, Vol.\ 66, pp.\ 15-37.

\ref[Bla86] Blackman, S.\ S., 1986.\ Multi-Target Tracking with Radar
Applications, Artech House, Dehdam, MA.

\ref[CBH09] Choi, H.\ L., Brunet, L., and How, J.\ P., 2009.\ ``Consensus-Based Decentralized Auctions for Robust Task Allocation," IEEE Transactions on Robotics, Vol.\ 25, pp.\ 912-926.

\ref[CGC04] Chang, H.\ S., Givan, R.\ L., and Chong, E.\ K.\ P., 2004.\ ``Parallel
Rollout for Online Solution of Partially Observable Markov Decision Processes," 
Discrete Event Dynamic
Systems, Vol.\ 14, pp.\ 309-341.

\ref[CXL19] Chu, Z., Xu, Z., and Li, H., 2019.\ ``New Heuristics for the RCPSP with Multiple Overlapping Modes," Computers and Industrial Engineering, Vol.\ 131, pp.\ 146-156.

\ref [Cas93]
Casta\~non, D.\ A., 1993.\ ``Reverse Auction Algorithms for Assignment Problems," in
Algorithms for Network Flows and Matching, Johnson, D.\ S., and McGeoch, C.\ C.\ (eds.),
American Math.\ Soc., Providence, RI,  pp.\ 407-429.

\ref [Chr97] Christodouleas, J.\ D., 1997.\ ``Solution Methods for Multiprocessor Network Scheduling
Problems with Application to Railroad Operations," Ph.D.\ Thesis, Operations Research Center,
Massachusetts Institute of Technology. 

\ref[DuV99] Duin, C., and Voss, S., 1999.\ ``The Pilot Method: A Strategy for Heuristic Repetition with Application to the Steiner Problem in Graphs," Networks: An International Journal, Vol.\ 34, pp.\ 181-191.

\ref[FeV02] Ferris, M.\ C., and Voelker, M. M., 2002.\ ``Neuro-Dynamic Programming for
Radiation Treatment Planning," Numerical Analysis Group Research Report NA-02/06, 
Oxford University Computing Laboratory, Oxford University. 

\ref[FeV04] Ferris, M.\ C., and Voelker, M. M., 2004.\ ``Fractionation in Radiation
Treatment Planning," Mathematical Programming B, Vol.\ 102, pp.\ 387-413.

\ref[GDM19] Guerriero, F., Di Puglia Pugliese, L., and Macrina, G., 2019.\ ``A Rollout Algorithm for the Resource Constrained Elementary Shortest Path Problem," Optimization Methods and Software, Vol.\ 34, pp.\ 1056-1074.

\ref[GTO15] Goodson, J.\ C., Thomas, B.\ W., and Ohlmann, J.\ W., 2015.\  ``Restocking-Based Rollout Policies for the Vehicle Routing Problem with Stochastic Demand and Duration Limits," Transportation Science, Vol.\ 50, pp.\ 591-607.

\ref[GuM01] Guerriero, F., and Musmanno, R., 2001.\ ``Label Correcting Methods to
Solve Multicriteria Shortest Path Problems," J.\ Optimization Theory Appl., Vol.\ 111,
pp.\ 589-613.

\ref[GuM03] Guerriero, F., and Mancini, M., 2003.\ ``A Cooperative Parallel Rollout
Algorithm for the Sequential Ordering Problem," Parallel Computing, Vol.\ 29, pp.\ 663-677.

\ref[HJG16] Huang, Q., Jia, Q.\ S., and Guan, X., 2016.\ ``Robust Scheduling of EV Charging Load with Uncertain Wind Power Integration," IEEE Trans.\  on Smart Grid, Vol.\ 9, pp.\ 1043-1054.

\ref[Jaf84] Jaffe, J.\ M., 1984.\ ``Algorithms for Finding Paths with Multiple
Constraints," Networks, Vol.\ 14, pp.\ 95-116.

\ref[KAH15] Khashooei, B.\ A., Antunes, D.\ J. and Heemels, W.P.M.H., 2015.\ ``Rollout Strategies for Output-Based Event-Triggered Control," in Proc.\ 2015 European Control Conference, pp.\ 2168-2173.

\ref[LGW16] Lan, Y., Guan, X., and Wu, J., 2016.\ ``Rollout Strategies for Real-Time Multi-Energy Scheduling in Microgrid with Storage System," IET Generation, Transmission and Distribution, Vol.\ 10, pp.\ 688-696.

\ref[LiW15] Li, H., and Womer, N.\ K., 2015.\ ``Solving Stochastic Resource-Constrained Project Scheduling Problems by Closed-Loop Approximate Dynamic Programming.," European J.\ of Operational Research, Vol.\ 246, pp.\ 20-33.

\ref [MMB02] McGovern, A., Moss, E., and Barto, A., 2002.\ ``Building a Basic Building
Block Scheduler Using Reinforcement Learning and Rollouts,"  Machine Learning, Vol.\
49, pp.\ 141-160.

\ref[MaJ15] Mastin, A., and Jaillet, P., 2015.\ ``Average-Case Performance of Rollout Algorithms for Knapsack Problems," J.\  of Optimization Theory and Applications, Vol.\ 165, pp.\ 964-984.

\ref [Mar84] Martins, E.\ Q.\ V., 1984.\ ``On a Multicriteria Shortest Path Problem,"
European J.\ of Operational Research, Vol.\ 16, pp.\ 236-245.

\ref [MPP04] Meloni, C., Pacciarelli, D., and Pranzo, M., 2004.\ 
``A Rollout Metaheuristic
for Job Shop Scheduling Problems," Annals of Operations Research, Vol.\
131, pp.\ 215-235.

\ref [PBW92] Pattipati, K.\ R., Deb, S.,  Bar-Shalom, Y., and Washburn, R.\ B., 1992.\ ``A New
Relaxation Algorithm and Passive Sensor Data Association,"  IEEE Trans.\ Automatic
Control, Vol.\ 37, pp.\ 198-213.

\ref[PPB01] Popp, R.\ L., Pattipati, K.\ R., and Bar-Shalom, Y., 2001.\ ``$m$-Best SD Assignment Algorithm with Application to Multitarget Tracking," IEEE Transactions on Aerospace and Electronic Systems, Vol.\ 37, pp.\ 22-39.

\ref [PoR97] Poore, A.\ B., and Robertson, A.\ J.\ A., 1997.\ New Lagrangian Relaxation Based Algorithm for a
Class of Multidimensional Assignment Problems," Computational Optimization and Applications, 
Vol.\ 8, pp.\ 129-150.

\ref [Poo94] Poore, A.\ B., 1994.\ ``Multidimensional Assignment Formulation of Data Association Problems Arising
from Multitarget Tracking and Multisensor Data Fusion," Computational Optimization and Applications, Vol.\
3, pp.\ 27-57.

\ref[SHB15] Simroth, A., Holfeld, D., and Brunsch, R., 2015.\ ``Job Shop Production Planning under Uncertainty: A Monte Carlo Rollout Approach," Proc.\ of the International Scientific and Practical Conference, Vol.\ 3, pp. 175-179.

\ref [SGC02] Savagaonkar, U., Givan, R., and Chong, E.\ K.\ P., 2002.\ ``Sampling
Techniques for Zero-Sum, Discounted Markov Games," in Proc.\ 40th 
Allerton Conference on Communication, Control and Computing, Monticello, Ill.

\ref[SLJ13] Sun, B., Luh, P.\ B., Jia, Q.\ S., Jiang, Z., Wang, F., and Song, C., 2013.\ ``Building Energy Management: Integrated Control of Active and Passive Heating, Cooling, Lighting, Shading, and Ventilation Systems," IEEE Trans.\ on Automation Science and Engineering, Vol.\ 10, pp.\ 588-602.

\ref[SNC18] Sarkale, Y., Nozhati, S., Chong, E.\ K., Ellingwood, B.\ R., and Mahmoud, H., 2018.\ ``Solving Markov Decision Processes for Network-Level Post-Hazard Recovery via Simulation Optimization and Rollout," in 2018 IEEE 14th International Conference on Automation Science and Engineering, pp.\ 906-912.

\ref[SZL08] Sun, T., Zhao, Q., Lun, P., and Tomastik, R., 2008.\ ``Optimization of Joint Replacement Policies for Multipart Systems by a Rollout Framework," IEEE Trans.\  on Automation Science and Engineering, Vol.\ 5, pp.\ 609-619.

\ref [Sec00] Secomandi, N., 2000.\ ``Comparing Neuro-Dynamic Programming Algorithms for
the Vehicle Routing Problem with Stochastic Demands,"
Computers and Operations Research, Vol.\ 27, pp.\ 1201-1225.

\ref [Sec01] Secomandi, N., 2001.\ ``A Rollout Policy for the Vehicle Routing Problem
with Stochastic Demands," Operations Research, Vol.\ 49, pp.\ 796-802.

\ref [Sec03] Secomandi, N., 2003.\ ``Analysis of a Rollout Approach to Sequencing
Problems with Stochastic Routing Applications," J.\ of Heuristics, Vol.\ 9, pp.\ 321-352.

\ref[StW91] Stewart, B.\ S., and  
White, C.\ C., 1991.\ ``Multiobjective $A^*$," J.\ ACM, Vol.\ 38, pp.\ 775-814.Ê

\ref[TGL13] Tesauro, G., Gondek, D.\ C., Lenchner, J., Fan, J., and Prager, J.\ M., 2013.\ ``Analysis of Watson's Strategies for Playing Jeopardy!," J.\  of Artificial Intelligence Research, Vol.\ 47, pp.\ 205-251.

\ref [TPL04] Tu, F.,  Pham, D., Luo, J., Pattipati, K.\ R., and Willett, P., 2004.\ 
``Decision Feedback with
Rollout for Multiuser Detection in Synchronous CDMA," IEE Proceedings-Communications, 
 Vol.\ 151, pp.\ 383-386.

\ref[TeG96] Tesauro, G., and Galperin, G.\ R., 1996.\ ``On-Line Policy Improvement
Using Monte Carlo Search,'' presented at the 1996 Neural Information
Processing Systems Conference, Denver, CO; also in M. Mozer et al.\ (eds.), Advances in
Neural Information Processing Systems 9, MIT Press (1997).

\ref[TsB00] Tseng, P., and Bertsekas, D.\ P., 2000.\ ``An $\e$-Relaxation Method for Separable Convex Cost Generalized Network Flow Problems," Math.\ Programming, Vol.\ 88, pp.\ 85-104.

\ref [TuP03] Tu, F., and Pattipati, K.\ R., 2003.\ ``Rollout Strategies for Sequential
Fault Diagnosis," IEEE Trans.\ on Systems, Man and Cybernetics, Part A, 
pp.\  86-99.

\ref[UGM18] Ulmer, M.\ W., Goodson, J.\ C., Mattfeld, D.\ C., and Hennig, M., 2018.\ ``Offline-Online Approximate Dynamic Programming for Dynamic Vehicle Routing with Stochastic Requests," Transportation Science, Vol.\ 53, pp.\ 185-202.

\ref[Ulm17] Ulmer, M.\ W., 2017.\ Approximate Dynamic Programming for Dynamic Vehicle Routing, Springer, Berlin.

\ref[WCG03] Wu, G.,  Chong, E.\ K.\ P., and Givan, R.\ L., 2003.\ ``Congestion Control
Using Policy Rollout," Proc.\ 2nd IEEE CDC, Hawaii, pp.\
4825-4830.

\ref[YDR04] Yan, X., Diaconis, P., Rusmevichientong, P., and Van Roy, B., 2004.\
``Solitaire: Man Versus Machine,'' Advances in Neural Information
Processing Systems, Vol.\ 17, pp.\ 1553-1560.


\vfill\eject

\pn{\bf APPENDIX: An Overview of the Auction Algorithm for 2-Dimensional Assignment}
\smskip


\pn 	In this appendix, we describe briefly the auction algorithm, an
intuitive method for solving the classical 2-dimensional assignment problem, which outperforms substantially its main competitors, both in theory and in practice, and it
can naturally exploit the availability of parallel computation. It is particularly well-suited for contexts of frequent reoptimization, where many problems are solved with slightly different data, as in the rollout algorithm of the present paper when applied to multidimensional assignment. The reason is that the final solution obtained for a given assignment problem by the auction algorithm (a set of prices, as we will describe shortly) is a good starting point for applying the algorithm to a similar problem.

In this appendix, we will sketch the basic principles of the auction algorithm, we will explain its
computational properties, and we will discuss some of its extensions to more general network flow
problems. These extensions are well-suited for use in rollout algorithms for complex combinatorial optimization problems beyond the multidimensional assignment context. For detailed presentations, we refer to 
the author's textbooks [Ber91], [Ber98], and the survey paper [Ber92]. For an extensive computational study, we refer
to Casta\~non [Cas93]. Several coded implementations of the auction algorithm are freely available from the internet, including some from the author's web site. The algorithm was first proposed in a 1979 report by the author
[Ber79].  The present appendix is based on the author's brief tutorial article [Ber01].

In the
classical assignment problem there are $n$ persons and $n$ objects that we have
to match on a one-to-one basis. There is  a benefit $a_{ij}$ for matching person
$i$ with object $j$ and we want to assign persons to objects so as to maximize
the total benefit.\footnote{\dag}{\ninepoint The auction algorithm is more intuitively explained in a context where we want to maximize benefit rather than minimizing cost. We will thus adopt the maximization context, but we can simply convert to minimization by changing the sign of $a_{ij}$. } Mathematically, we want to find a one-to-one assignment [a set of person-object pairs $(1,j_1),\ldots,(n,j_n)$, such that the objects
$j_1,\ldots,j_n$ are all distinct] that maximizes the total
benefit $\sum_{i=1}^na_{ij_i}$.

The assignment problem is important in many practical contexts. The  most
obvious ones are resource allocation problems, such as assigning personnel to
jobs, machines to tasks, and the like. There are also many situations where the
assignment problem appears as a subproblem in various methods for solving more
complex problems, as in the context of the present paper. 

The assignment problem is also of great theoretical importance because, despite
its simplicity, it embodies a fundamental network optimization structure.  All linear single commodity network flow problems,
can be reduced to the assignment problem by means of a simple reformulation.
Thus, any method for solving the assignment problem can be generalized to solve
the linear network flow problem, and in fact this approach is particularly helpful in
understanding the extension of auction algorithms to network flow problems that are more
general than assignment.

The principles on which the auction algorithm is based differ markedly from all other linear network flow methods. In particular, classical methods for assignment are based on iterative improvement of some cost
function; for example a primal cost (as in primal simplex methods), or a dual cost (as in
Hungarian-like methods,  dual simplex methods, and relaxation methods). The auction
algorithm departs significantly from the cost improvement idea; at any one iteration,  it
may deteriorate both the primal and the dual cost, although in the end it finds an optimal
assignment. It is based on a notion of approximate optimality, called {\it
$\e$-complementary slackness\/}, and while it implicitly tries to solve a dual problem, it
actually attains a dual solution that is not quite optimal. 

\subsection{The Auction Process}

\pn To develop
an intuitive understanding of the auction algorithm, it is helpful to introduce
an economic equilibrium problem that turns out to be equivalent to the
assignment problem.
Let us consider the possibility of matching the $n$ objects with the  $n$ persons
through a market mechanism, viewing each person as an economic agent acting in
his own best interest. Suppose that object $j$ has a price $p_j$ and that the
person who receives the object must pay the price $p_j$. Then, the (net) value
of object $j$ for person $i$ is $a_{ij}-p_j$ and each person $i$ would logically
want to be assigned to an object $j_i$ with maximal value, that is, with
$$a_{ij_i}-p_{j_i}=\max_{j=1,\ldots,n}\{a_{ij}-p_j\}.\eqno(1)$$ 
We will say that
a person $i$ is {\it happy}  if this condition holds  and we will say that an
assignment and a set of prices are at {\it equilibrium} when all persons are
happy.

Equilibrium assignments and prices are naturally of great interest to
econo\-mists, but there is also a fundamental relation with the assignment
problem; it turns out that an equilibrium assignment offers maximum total
benefit (and thus solves the assignment
problem), while the corresponding set of prices solves an associated dual
optimization problem. This is a consequence  of the celebrated duality theorem
of linear programming.

Let us consider now a natural process for finding an
equilibrium assignment. We will call this process the {\it naive auction
algorithm\/}, because it has a serious flaw, as will  be seen shortly.
Nonetheless, this flaw will help motivate a more sophisticated and correct
algorithm.

The naive auction algorithm proceeds in ``rounds'' (or ``iterations'') starting
with {\it any} assignment and {\it any} set of prices. There is an assignment and
a set of prices at the beginning of each
round, and if all persons are happy with these, the process terminates.
Otherwise some person who is not happy is selected. This person, call him $i$,
finds an object $j_i$ which offers maximal value, that is,
$$j_i\in\arg\max_{j=1,\ldots,n}\{a_{ij}-p_j\},\eqno(2)$$ 
and then:

\nitem{(1)} Exchanges
objects with the person assigned to $j_i$ at the beginning of the round.

\nitem{(2)} Sets the price of the best object $j_i$ to the
level at which he is
indifferent between $j_i$ and the second best object, that is, he sets $p_{j_i}$ to 
$$p_{j_i}+\g_i,\eqno(3)$$
where 
$$\g_i=v_i-w_i,\qquad \eqno(4)$$ 
$v_i$ is the best object
value, $$v_i=\max_j\{a_{ij}-p_j\},\qquad \eqno(5)$$
and $w_i$ is the second best object value
$$w_i=\max_{j\ne j_i}\{a_{ij}-p_j\}, \qquad \eqno(6)$$
that is, the best value over objects other than $j_i$. (Note that $\g_i$ is the
largest increment by which the best object price $p_{j_i}$ can be increased, with
$j_i$ still being the best object for person $i$.)
\smskip
 
\pn This process is repeated in a
sequence of rounds until all persons are happy.

We may view this process as an auction, where at each round the bidder $i$
raises the price of his or her preferred object by the {\it bidding
increment} $\g_i$.  Note that $\g_i$
cannot be negative since $v_i\ge w_i$ [compare Eqs.\ (5) and (6)], so the
object prices tend to increase.
Just as in a real auction, bidding increments and price increases spur
competition by making the bidder's own preferred object less attractive to other
potential bidders.

Unfortunately, this auction does not always  work. The difficulty
is that the bidding increment $\g_i$ is zero when more than one
object offers maximum value for the bidder $i$ [cf.\ Eqs.\ (4), (6)]. As a result, a
situation may be created where several persons contest a smaller number of equally
desirable objects without raising their prices, thereby creating a never ending cycle.

To break such cycles, we introduce a perturbation mechanism, motivated by real
auctions where each bid for an object must raise its price by a minimum
positive increment, and bidders must on occasion take risks to win their
preferred objects. In particular, let us fix a positive scalar $\e$ and say that
a person $i$ is {\it almost happy}  with an assignment and a set of prices if
the value of its assigned object $j_i$ is within $\e$ of being maximal, that is,
$$a_{ij_i}-p_{j_i}\ge\max_{j=1,\ldots,n}\{a_{ij}-p_j\}-\e.\eqno(7)$$  
We will
say that an assignment and a set of prices are {\it almost at equilibrium} when
all persons are almost happy. The condition (7), introduced first in 1979
in conjunction with the auction algorithm, is known as {\it $\e$-complementary
slackness\/} and plays a central role in several optimization contexts. For $\e=0$ it
reduces to ordinary complementary slackness [compare Eq.\ (1)].

We now reformulate the previous auction process so that the bidding
increment is always at least equal to $\e$. The resulting method, the {\it
auction algorithm\/}, is the same as the naive auction algorithm, except that the
bidding increment $\g_i$ is
$$\g_i=v_i-w_i+\e,\qquad \eqno(8)$$
[rather than $\g_i=v_i-w_i$ as in Eq.\ (4)].
With this choice, the bidder of a round is almost happy at the end of the round
(rather than happy). The particular increment
$\g_i=v_i-w_i+\e$ used in the auction algorithm is the maximum amount with this
property. Smaller increments $\g_i$ would also work as long as $\g_i\ge\e$, but
using the largest possible increment accelerates the algorithm. This is
consistent with experience from real auctions, which tend to terminate faster
when the bidding is aggressive.

 We can now show that this reformulated auction process terminates in a finite
number of rounds, necessarily with an assignment and a set of prices that are
almost at equilibrium. To see this, note that once an object receives a bid
for the first time, then the person assigned to the object at every subsequent
round is almost happy; the reason is that a person is almost happy just after
acquiring an object through a bid, and continues to be almost happy as long as he
holds the object (since the other object prices cannot decrease in the course of
the algorithm). Therefore, the persons that are not almost happy must be assigned to
objects that have never received a bid. In particular, {\it once each object
receives at least one bid, the algorithm must terminate\/}. Next note that  if an
object receives a bid in $m$ rounds, its price must exceed its initial price by at
least $m\e$. Thus, for sufficiently large $m$, the object will become ``expensive''
enough to be judged ``inferior'' to some object that has not received a bid so far.
It follows that only for a limited number of rounds can an object receive
a bid while some other object still has not yet received any bid.
Therefore, there are two possibilities: either (a) the auction terminates in a
finite number of rounds, with all persons almost happy, before every object receives
a bid or (b) the auction continues until, after a finite number of rounds, all
objects receive at least one bid, at which time the auction terminates. (This
argument assumes that any person can bid for any object, but it can be
generalized for the case where the set of feasible person-object pairs is
limited, as long as at least one feasible assignment exists.)

\vskip-0.1pc

\subsection{\bf Optimality Properties at Termination}\medskip
\vskip-0.5pc

\pn 
When the auction algorithm terminates, we have an assignment that
is almost at equilibrium, but does this assignment maximize the total benefit?
The answer here depends strongly on the size of $\e$. In a real auction, a
prudent bidder would not place an excessively high bid for fear that he might
win the object at an unnecessarily high price. Consistent with this intuition, we
can show that if $\e$ is small, then the final assignment will be ``almost
optimal.'' In particular, we can show that  {\it the total benefit of the final
assignment is within $n\e$ of being optimal\/}. To see this, note that an
assignment and a set of prices that are almost at equilibrium may be viewed as being at
equilibrium for a {\it slightly different} problem where all benefits $a_{ij}$ are the
same as before, except for the $n$ benefits of the assigned pairs which are modified by an
amount no more than $\e$.

Suppose now that the benefits $a_{ij}$ are all integer, which is the typical
practical case (if $a_{ij}$ are rational numbers, they can be scaled up to
integer by multiplication with a suitable common number). Then, the
total benefit of any assignment is integer, so if $n\e<1$, a complete assignment
that is within $n\e$ of being optimal must be optimal. It follows, that {\it if
$$\e<{1\over
n},$$
and the benefits
$a_{ij}$ are all integer, then the assignment obtained
upon termination of the auction algorithm is optimal\/}. Let us also note that the
final set of prices is within $n\e$ of being an optimal solution of the dual
problem
$$\min_{p_j\atop j=1,\ldots
,n}\lf\{\sum_{j=1}^np_j+\sum_{i=1}^n\max_j\bl\{a_{ij}-p_j\br\}\ri\}.\eqno(9)$$
This leads to the interpretation of the auction algorithm as a dual algorithm (in fact an
approximate coordinate ascent algorithm; see the cited literature).

\subsection{\bf $\e$-Scaling and Reoptimization}\medskip

\pn The auction algorithm exhibits interesting computational behavior, and it is
essential to understand this behavior in order to use and implement the algorithm
efficiently.
First note that the amount of work to solve the problem can depend strongly
on the value of $\e$ and on the maximum absolute object value
$$C=\max_{i,j}|a_{ij}|.$$
Basically, for many types of problems, the number of bidding rounds up to
termination tends to be proportional to $C/\e$. Note also that there is a dependence on the
initial prices; if these prices are ``near optimal,'' we expect that the number of rounds
to solve the problem will be relatively small.

The preceding observations suggest the idea of {\it $\e$-scaling}, proposed in the original paper [Ber79]], which
consists of applying the algorithm several times, starting with a large value of
$\e$ and successively reducing $\e$ up to an ultimate value that is less than
some critical value (for example, $1/n$, when the benefits $a_{ij}$ are integer).
Each application of the algorithm provides good initial prices for the next
application. This is a common idea in nonlinear programming, encountered for
example, in barrier and penalty function methods. An alternative form of scaling,
called {\it cost scaling\/}, is based on successively representing the benefits
$a_{ij}$ with an increasing number of bits, while keeping $\e$ at a constant
value. 

In practice,  it is a good idea to at least consider scaling. For sparse
assignment problems, that is, problems where the set of feasible assignment
pairs is severely restricted, scaling seems almost universally helpful. In theory, scaling
leads to auction algorithms with a particularly favorable polynomial complexity (without
scaling, the algorithm is pseudopolynomial; see the cited literature).

Finally, let us note that the auction algorithm can be initialized with arbitrary prices, as well arbitrary $\e$. Within this context, it is important to choose initial prices that are close to the final prices obtained upon termination. This can be very important in a reoptimization setting, and in fact this is born out from computational complexity analysis (see [Ber88] or [Ber98]). Indeed, reoptimization with a favorable set of starting prices underlies the idea of $\e$-scaling or cost scaling. The rollout framework of this paper involves extensive reoptimization, with solutions of slightly differing assignment problems. This suggests an overwhelming advantage of the auction algorithm as a base heuristic within the context of this paper.

\subsection{\bf Parallel and Asynchronous Implementation}\medskip

\pn Both the bidding and the assignment phases of the auction
algorithm are highly parallelizable. In particular, the bidding and the assignment
can be carried out for all persons and objects simultaneously.
Such an implementation can be termed {\it synchronous}. There are also {\it
totally asynchronous\/} implementations of the auction algorithm, which are
interesting because they are quite flexible and also tend to result in faster
solution in some types of parallel machines. To understand these
implementations, it is useful to think of a person as an autonomous decision
maker who at unpredictable times obtains information about the prices of the
objects. Each person who is not almost happy makes a bid at arbitrary times on
the basis of its current object price information (that may be outdated because
of communication delays). 

Bertsekas and Casta\~non [BeC91] give a careful
formulation of the totally asynchronous model, and a proof of its validity. They
include also extensive computational results on a shared memory machine,
confirming the advantage of asynchronous over synchronous implementations. In the context of rollout, there is potential for much more parallelization in the context of rollout, since as many as $m$ independent 2-dimensional assignment problems can be solved in parallel at each stage of the rollout computations.

\subsection{\bf Variations and Extensions}
\medskip

\pn The auction algorithm can be extended to solve a number of variations of the assignment
problem, such as the asymmetric assignment problem where
the number of objects is larger than the number of persons and there is a requirement that
all persons be assigned to some object. Naturally, the notion of an assignment must now be
modified appropriately. To solve this problem, the auction algorithm need only be modified
in the choice of initial conditions. It is sufficient to require that all initial prices be
zero. A similar algorithm can be used for the case where there is no requirement that all
persons be assigned.

There have been extensions of the auction algorithm for other types of linear network
optimization problems. These extensions find potential application in the context of rollout algorithms for challenging combinatorial problems, which involve graphs and network flows. The general approach for constructing  auction algorithms for such problems is to convert them to
assignment problems, and then to suitably apply the auction algorithm and
streamline the computations. In particular, the classical shortest path
problem can be solved correctly by the naive auction algorithm described earlier, once the
method is streamlined. Similarly, auction algorithms can be constructed for the max-flow
problems, and are very efficient. These algorithms bear a close relation to preflow-push
algorithms for the max-flow problem.

Linear transportation problems may also be addressed with the auction algorithm (Bertsekas and Casta\~non [Ber89]). The basic idea is to convert the transportation problem into an assignment problem by creating multiple copies of persons (or objects) for each source (or sink
respectively), and then to modify the auction algorithm to take advantage of the
presence of the multiple copies, while appropriately streamlining the computations. This auction algorithm may be used in the context of the facility location problem of Example 3.2.

There are also extensions of the auction algorithm for linear
minimum cost flow (transshipment) problems, such as the so called $\e$-relaxation
method, and the auction/sequential shortest path algorithm algorithm (see the network optimization textbook [Ber98] for a detailed description and further references).  The $\e$-relaxation method was first published by the author in [Ber86], although it was known much earlier (since the development of the
mathematically equivalent auction algorithm). It is equivalent to the so called pre-flow push algorithms, as discussed in the author's paper [Ber93]. These methods have interesting theoretical properties and like the auction algorithm, are well suited for parallelization (see the papers [Ber86], [BeG97], the survey by Bertsekas, Casta\~ non,  Eckstein, and Zenios [BCE95], and the textbook by Bertsekas and Tsitsiklis [BeT89]).

Let us finally note that there have been proposals of auction algorithms for convex
separable network optimization problems with and without gains (but with a single
commodity and without side constraints); see Tseng and Bertsekas [TsB00].

\end